\documentclass[12pt]{amsart}
\usepackage{amssymb,latexsym,amsthm,mathptmx}
\usepackage{amsmath,amsfonts,amssymb}
\usepackage{graphicx}
\usepackage{color}
\usepackage{mathrsfs}
\usepackage{eepic}
\usepackage{epic}
\newtheorem{thm}{Theorem}[section]
\newtheorem{crl}[thm]{Corollary}

\newtheorem{lmm}[thm]{Lemma}

\newtheorem{prp}[thm]{Proposition}

\newtheorem{thmak}{The Arens-Kelley theorem}

\theoremstyle{definition}
\newtheorem{dfn}[thm]{Definition}
\newtheorem{exa}[thm]{Example}
\newtheorem{newProof}{\it{Proof of Proposition \ref{alcreg}}}
\renewcommand{\thenewProof}

\theoremstyle{remark}
\newcommand{\bA}{\mathcal{A}}

\newcommand\R{\mathbb R}
\newcommand\C{\mathbb C}

\newcommand\D{\mathbb D}

\newcommand{\K}{\mathbb{K}}

\newcommand{\F}{\frak{F}}
\newcommand{\T}{\mathbb{T}}
\newcommand{\TT}{\mathfrak{T}}

\newcommand{\suppm}{\operatorname{supp}}

\newcommand{\hh}{\mathfrak{h}}

\newcommand{\ball}{\operatorname{Ball}}

\newcommand{\fraca}{\frac{\alpha}{|\alpha|}}

\newcommand{\ex}{\operatorname{ext}}

\newcommand{\ch}{\operatorname{Ch}}

\newcommand{\ext}{\operatorname{ext}}

\newcommand{\rpi}{\operatorname{Ran}_{\pi}}
\begin{document}
\title[The Mazur-Ulam property for a Banach space]
{The Mazur-Ulam property for a Banach space which satisfies a separation condition
}
\author{
Osamu~Hatori
}
\address{
Institute of Science and Technology,
Niigata University, Niigata 950-2181, Japan
}
\email{hatori@math.sc.niigata-u.ac.jp
}
%
%
%
%
%
%
%
%


\begin{abstract}      
We study $C$-rich spaces, lush spaces, and $C$-extremely regular spaces  concerning with the Mazur-Ulam property. We show that a uniform algebra and the real part of a uniform algebra with the supremum norm are  $C$-rich spaces, hence  lush spaces. We prove that a uniformly closed subalgebra of the algebra of complex-valued continuous functions on a locally compact Hausdorff space which vanish at infinity is $C$-extremely regular provided that it separates the points of the underlying space and has no common zeros. 
In section 3 we exhibit descriptions on the Choquet bounday, the \v Silov bounday, the strong boundary points. We also recall the definition that a function space strongly separates the points in the underlying space. We need to avoid the confusion which appears because of the variety of names of these concepts;  they  sometimes differs from authors to authors. 
After some preparation, 
we study the Mazur-Ulam property in sections \ref{sec4} through \ref{sec6}. We exhibit a sufficient condition on a Banach space which has the Mazur-Ulam property and the complex Mazur-Ulam property (Propositions \ref{prop2.2} and \ref{prop15}). 
In section \ref{sec5} we consider a Banach space with a separation condition $(*)$ (Definition \ref{dfn*}). We prove that a real Banach space satisfying $(*)$  has the Mazur-Ulam propety (Theorem \ref{real*}), and a complex Banach space satisfying $(*)$ has the complex Mazur-Ulam property (Theorem \ref{71}). 
Applying the results in the previous sections we prove that an extremely $C$-regular complex linear subspace has the complex Mazur-Ulam property  (Corollary \ref{ecregmu}) in section \ref{sec6}. As a consequence 
we prove that any closed subalgebra of the algebra of all complex-valued continuous functions defined on a locally compact Hausdorff space has the complex Mazur-Ulam property (Corollary \ref{csacmu}).
\end{abstract}

\maketitle

\section{Introduction}\label{sec1}
Tingley's problem asks whether every  surjective isometry between the unit spheres of Banach spaces is extended to a surjective isometry between whole spaces. Tingley \cite{tingley} raised the problem  in 1987.
First solution of Tingley's problem seems to be by 
Wang \cite{wang}.
who dealt with the space of all $\K$-valued ($\K=\R$ or $\C$) continuous functions which vanish at infinity on a locally compact Hausdorff space $Y$ (cf. \cite{wang1996a}).  
A considerable number of interesting  results have shown that Tingley's problem has an affirmative answer, and
no counterexample is known.
According to  \cite[p.730]{yangzhao}, Ding was the first to consider Tingley's problem between different type of spaces \cite{ding2003B}.  In fact,  Ding 
\cite[Corollary 2]{ding2007} proved that the real Banach space of all null sequences of real numbers satisfies what we now call the Mazur-Ulam property. Liu \cite{liu2007} also has an early contribution to the Tingley's problem on different types of spaces.  Later
Cheng and Dong \cite{chengdong} formally introduced the concept of the Mazur-Ulam property.

A real Banach space $B$ has the {\em Mazur-Ulam property} if for any real Banach space $B'$, every surjective isometry from the unit sphere of $B$ onto the unit sphere of $B'$ admits an extension to a surjective real-linear isometry from $B$ onto $B'$.
Tan \cite{tan2011a,tan2011b,tan2012a} showed that the space $L^p(\R)$ for $\sigma$-finite positive measure space has the Mazur-Ulam property. 
Boyko, Kadates, Mart\'in and Werner \cite[Theorem 2.4]{bkmw} introduced  $C$-richness and lushness for  subspaces of continuous functions and proved that a $C$-rich subspace is lush. 
Tan, Huang and Liu \cite{thl2013} introduced the notion of 
local GL (generalized lush) spaces and proved that every local GL space has the Mazur-Ulam property.

Tanaka \cite{tanaka} opened another direction in the study of Tingley's problem by exhibiting a positive solution  
for the Banach algebra of complex matrices. 
Mori and Ozawa \cite{moriozawa} proved that the Mazur-Ulam  property  for unital $C^*$-algebras and real von Neumann algebras. 
Cueto-Avellaneda and Peralta \cite{cp2019} proved that a complex (resp. real) Banach space of all continuous maps with the value in a complex (resp. real) Hilbert space has the Mazur-Ulam property (cf. \cite{cp2020}).
The results by 
Becerra-Guerrero, Cueto-Avellaneda, Fern\'andez-Polo and Peralta \cite{bcfp}  and 
Kalenda and Peralta \cite{kape} proved that any
 JBW*-triples have the Mazur-Ulam property. 

The Mazur-Ulam property for  a Banach space of dimension 2 has been unsolved for many years. 
The final solution was exhibited by the remarkable outstanding advance of Banakh \cite{banakh} who proved that any Banach space of dimension 2 has the Mazur-Ulam property. The problem on a Banach space of a finite dimension greater than 2 seems to be still open. 
The study
of the Mazur-Ulam property  is nowadays a
challenging subject of study (cf. \cite{banakh,cabellosanchez2019,cp2020,wh2019}). 

A complex Banach space $B$ is said to  have the {\em complex  Mazur-Ulam property}, emphasizing the term `complex',  
 if for any complex Banach space $B'$, every surjective isometry from the unit sphere of $B$ onto the unit sphere of $B'$ admits an extension to a surjective real-linear isometry from $B$ onto $B'$ (cf. \cite{hatori}).
Note that  a complex Banach space has the complex Mazur-Ulam property provided that it has the Mazur-Ulam property as a real Banach space  since a complex Banach space is a real Banach space simultaneously. 
Jim\'enez-Vargas, Morales-Compoy, Peralta and Ram\'irez \cite[Theorems 3.8, 3.9]{jmpr2019} probably  provided the first example of complex Banach spaces  having the complex Mazur-Ulam property (cf. \cite{peralta2019a}).

The complex Mazur-Ulam property for uniform algebras was proved in \cite{hatori}. 
The existence of unity in a uniform algebra is a key point for the proof of the property. 
The complex Mazur-Ulam property for a uniformly closed algebra on a locally compact Hausdorff space is a problem in \cite{hatori}.
Cueto-Avellaneda, Hirota, Miura and Peralta \cite{chmp} recently showed that each surjective
isometry between the unit spheres of two uniformly closed  algebras on locally compact
Hausdorff spaces which separates the points without common zeros admits an extension to a surjective real linear isometry between these algebras.
Very recently, Cabezas, Cueto-Avellaneda, Hirota, Miura and Peralta \cite{cchmp}  proved the complex Mazur-Ulam property for  a commutative JB*-triple. Both results concerns the spaces of continuous functions without constants. Note that the Mazur-Ulam property, or even the complex Mazur-Ulam property for non-unital $C^*$-algebras seems to be still missing.

In this paper we further study the problem on the complex Mazur-Ulam property. 
We introduce a separation condition which is named $(*)$ for a Banach space in section \ref{sec5}. We prove that a real (resp. complex) Banach space which satisfies the  condition $(*)$ has the (resp. complex) Mazur-Ulam property. 
An extremely $C$-regular subspace satisfies the condition. It was introduced by Fleming and Jamison \cite[Definition 2.3.9]{fj1}, which is a generalization of an extremely regular subspace coined by Cengiz \cite{cengiz}. We prove that an extremely $C$-regular complex linear subspace has the complex Mazur-Ulam property (Corollary \ref{ecregmu}). As a consequence the complex Mazur-Ulam property for 
 a closed subalgebra of the algebra of all complex-valued continuous functions defined on a locally compact Hausdorff space is established  (Corollary \ref{csacmu}). 
 
We recall the notions of the Choquet boundary, the strong boundary, the strong separation of the points in the underlying space in section \ref{sec3}. We expect that many of results in section \ref{sec3} are a part of folklore, but for the sake of self-contained exposition and for the convenience of the readers who are not familiar with these concepts for function spaces without constants, we include complete proofs of results as possible as we can.

For a real or complex Banach space $B$ 
the unit sphere $\{a\in B:\|a\|=1\}$ of $B$ is denoted by $S(B)$ and the unit ball $\{a\in B: \|a\|\le 1\}$ by $\operatorname{Ball}(B)$. The set of all maximal convex subsets of $S(B)$ is denoted by $\F_B$.
We denote $\K=\R$ or $\C$. 
$\T=\{z\in {\mathbb K}:|z|=1\}$.
Throughout the paper $Y$ denotes a locally compact Hausdorff space and $X$ a compact Hausdorff space. 
The space of all $\K$-valued continuous function on $Y$ which vanish at infinity is denoted by $C_0(Y,\K)$. If $Y$ is compact, then we simply denote $C(Y,\K)$ instead of $C_0(Y,\K)$
The supremum norm is denoted by $\|\cdot\|_\infty$.
%


\section{$C$-rich spaces, lush spaces and extremely $C$-regular spaces}\label{sec X}

\subsection{$C$-richness, lushness, the numerical index and the Mazur-Ulam property.}
A $C$-rich subspace was introduced by Boyko, Kadets, Mart\'in and Werner \cite{bkmw}.
\begin{dfn}[\cite{bkmw}]
A closed $\K$-linear subspace $E$ of $C(X,\K)$, $\K=\C$ or $\R$,  is called {\em $C$-rich} if for every nonempty open subset $U$ of $X$ and $\varepsilon>0$, there exists a positive function $h_\varepsilon$ of norm $1$ with support inside $U$ such that the distance from $h_\varepsilon$ to $E$ is less than $\varepsilon$.
\end{dfn}
Suppose that $X$ is a compact Hausdorff space without isolated points and $p_1,\dots, p_n\in C(X,\K)^*$. 
Then $E=\bigcap_{j=1^n}p_j^{-1}(0)$ is a $C$-rich subspace of $C(X,\K)$ \cite[Proposition 2.5]{bkmw}.
Furthermore, if $X$ is perfect, then every  subspace of $C(X,\K)$ of codimension finite is $C$-rich since in this case $C$-richness is equivalent to richness given  \cite[Proposition 1.2]{kp}.
Another example is a uniform algebra.
We say that $A$ is a {\em uniform algebra} on a compact Hausdorff space $X$ if  $A$ is a closed subalgebra of $C(X,\C)$ which contains constants and separates the points of $X$ (see \cite{br,gamegame,stout}, where a uniform algebra is called a function algebra in \cite{br}).
\begin{prp}\label{uacr}
Let $A$ be a uniform algebra on a compact Hausdorff space $X$ and $S$ the \v Silov bounday\footnote{The definition for uniform algebras, see \cite{br}, or a general definition is given in Definition \ref{shilov}.}  for $A$.
Then $A|S$ is $C$-rich $\C$-subspace of $C(S,\C)$.
\end{prp}
\begin{proof}
Let $U$ be an open subset of  $S$ and $\varepsilon_0>0$ arbitrary. We may suppose that $U$ is a proper subset of $S$. Let $0<\varepsilon<\min\{1/2,\varepsilon_0/(\sqrt{5}+1)\}$.  It is known that the Choquet boundary $\ch(A)$ for $A$ is dense in $S$ and each point $x\in \ch(A)$ is a strong boundary point  \cite{br} (cf. Theorem \ref{chalg} and Corollary \ref{431}). 
Hence there exist $p\in U\cap\ch(A)$ and  $f\in A|S$ such that $f(p)=1=\|f\|_\infty$ and $|f|<1$ on $S\setminus U$. Since $A$ is closed under the multiplication, we may suppose that $|f|<1/2$ on $S\setminus U$. Put 
\[
\Delta=\{z\in \bar{D}: \operatorname{Re}z\ge0, \,\,|\operatorname{Im}z|\le\varepsilon\}.
\]
Then by the well known Carath\'eodory theorem (cf. \cite{pome}) there is a homeomorphism $\pi_\varepsilon:\bar{D}\to \Delta$ such that $\pi_\varepsilon$ is analytic from $D$ onto the interior of  $\Delta$. We may assume that  $\pi_\varepsilon(1)=1$ and $\pi_\varepsilon(\{z\in \bar{D}:|z|<1/2\})\subset \{z\in \Delta:\operatorname{Re}z<\varepsilon\}$. As $\pi_\varepsilon$ is uniformly approximated by analytic polynomials on $\bar{D}$, we have $g=\pi_\varepsilon\circ f\in A|S$. Note that $0\le\operatorname{Re}g\le 1$ on $S$.  By  Urysohn's lemma there exists a continuous function  $h:S\to [0,1]$ such that 
\begin{equation*}
h(y)=
\begin{cases}
0, & \text{if $\operatorname{Re}g(y)\le \varepsilon$}, 
\\
1,& \text{if $\operatorname{Re}g(y)\ge 2\varepsilon$}.
\end{cases}
\end{equation*}
If $y\in S\setminus U$, then $|f(y)|<1/2$, so $\operatorname{Re}g(y)=\operatorname{Re}(\pi_\varepsilon\circ f)(y)<\varepsilon$. Thus $hg=0$ on $S\setminus U$. 
As $g(p)=\pi_\varepsilon\circ f(p)=1$ we have that $hg(p)=1$. As $\|g\|_\infty=g(p)=1=h(p)=\|h\|_\infty$ we have $\|hg\|_\infty=1$. We show that $\|hg-g\|\le \sqrt{5}\varepsilon$. Let $y\in S$. If $\operatorname{Re}g(y)\ge 2\varepsilon$, then $h(y)=1$. Hence $(hg-g)(y)=0$. Suppose that $\operatorname{Re}g(y)\le 2\varepsilon$.  As $0\le h(y)\le 1$, $|h(y)-1|\le 1$. As $g(y)\in \Delta$ and $\operatorname{Re}g(y)\le 2\varepsilon$ we infer that $|g(y)|\le \sqrt{5}\varepsilon$. Hence 
\[
|(hg-g)(y)|=|g(y)||h(y)-1|\le \sqrt{5}\varepsilon.
\]
As $\operatorname{Re}g\ge 0$ on $X$, we have $0\le h\operatorname{Re}g\le 1$. Put $h_0=h\operatorname{Re}g$. Then $h_0$ is a positive function of norm $1$ since $h_0(p)=1$. As $h=0$ on $S\setminus U$, we have $h_0=0$ on $S\setminus U$. Thus the support of $h_0$ is inside of $U$. 
We have
\[
\|h_0-g\|_\infty\le \|hg-g\|_\infty+\|h\|_\infty\|\operatorname{Re}g-g\|_\infty\le
\sqrt{5}\varepsilon+\varepsilon<\varepsilon_0.
\]
Thus $d(h_0, A|S)<\varepsilon_0$
\end{proof}
In the same way we see the following. 
\begin{prp}\label{ruarr}
Suppose that $A$ is a uniform algebra on a  compact Hausdorff space $X$ and $S$ its \v Silov boundary. 
Let $E$ be the uniform closure $\overline{\operatorname{Re}A|S}$ of the real part of $A|S$. 
Then $E$ is a $C$-rich $\R$-subspace of $C(X,\R)$
\end{prp}
\begin{proof}
Let $U$ be an open subset $U$ of $S$ and $\varepsilon_0>0$ arbitrary. 
In fact, a given $U$ and $\varepsilon_0>0$, $h_0$ and $g\in A$ are  the same functions as in the proof of Proposition \ref{uacr} we have
$\|h_0-\operatorname{Re}g\|_\infty\le \varepsilon_0$.
\end{proof}

A lush space was introduced by Boyko, Kadets, Mart\'in and Werner \cite{bkmw}.
\begin{dfn}\label{lush}
Let $B$ be a $\K$-Banach space, where $\K=\R$ or $\C$.
Let $\delta>0$ and $p\in S(B^*)$. The {\em slice} denoted by $SL(\ball(B),p,\delta)$ is
\[
\{a\in \ball(B): \operatorname{Re}p(a)>1-\delta\}.
\]
$B$ is said to be {\em $\K$-lush} if for every $a,b\in S(B)$ and $\varepsilon>0$, there exists $q\in S(B^*)$ such that  $b\in S(\ball(B), q, \varepsilon)$ and 
\[
d(a, \operatorname{co}(\T SL(\ball(B), q, \varepsilon)))<\varepsilon,
\]
where $\T$ stands $\{\pm 1\}$ if $\K=\R$ and  $\{z\in \C:|z|=1\}$ if $\K=\C$, $\operatorname{co}(\cdot)$ stands the convex hull.
\end{dfn}
Boyko, Kadets, Mart\'in and Werner \cite[Theorem 2.4]{bkmw} proved that a $C$-rich $\K$-subspace of $C(X,\K)$ for a compact Hausdorff space $X$ is $\K$-lush. 
Tan, Huang and Liu introduced  a local-Gl-space in \cite{thl2013}, which is a real Banach space. 
They proved that $\R$-lush space is a local-GL-space \cite[Theorem 3.7]{thl2013} and every local-GL-space has the Mazur-Ulam property. 
Let us briefly recall that a real Banach space  $B$ is GL-{\em space} if for every $a\in S(B)$ and every $0<\varepsilon <1$ there exists a 
 $p\in S(B^*)$ such that 
\[
d(b,SL(\ball (B), p,\varepsilon))+d(-b,SL(\ball (B), p,\varepsilon))<2+\varepsilon
\]
for all $b\in S(B)$. A real Banach space $B$ is said to be a {\em local {\rm{GL}}-space} if for every separable subspace $E$ of $B$ , there exists a GL-subspace $E'$ such that $E\subset E'\subset B$. 
Note that an $\R$-lush space is a local GL-space \cite[Example 3.6]{thl2013}.
\begin{crl}
Every uniform algebra is $\C$-lush.
The uniform closure of the real part of a uniform algebra is $\R$-lush. Hence the uniform closure of the real part of a uniform algebra has the Mazur-Ulam property.
\end{crl}
\begin{proof}
Let $A$ be a uniform algebra on $X$ and $S$ the \v Silov boundary for $A$.
Then by Proposition \ref{uacr}, $A|S$ is $C$-rich $\C$-subspace of $C(S,\C)$.
Then by \cite[Theorem 2.4]{bkmw} we see that $A|S$ is $\C$-lush. Note that lushness is invariant under the isometries due to the definition of lushness. Therefore $A$ is $\C$-lush. 

By Proposition \ref{ruarr}  $\overline{\operatorname{Re}A|S}$ is $C$-rich. 
Then \cite[Theorem 2.4]{bkmw} ensures that $\overline{\operatorname{Re}A|S}$ is $\R$-lush. 
As the \v Silov boundary is a boundary, $\overline{\operatorname{Re}A|S}=\overline{\operatorname{Re}A}|S$. Since lushness is invariant under the isometries, we see that $\overline{\operatorname{Re}A}$ is $\R$-lush. Hence $\overline{\operatorname{Re}A}$ has the Mazur-Ulam property.
\end{proof}
It seems not to be known if a $\C$-lush space has the complex Mazur-Ulam property or not.
We proved that a uniform algebra has the complex Mazur-Ulam property in \cite{HOST}.
Recall that a uniform algebra $A$ on $X$ is a Dirichlet algebra provided that  $\overline{\operatorname{Re}A}=C(X,\R)$. Several uniform algebras including the disk algebra on the unit circle  is a Dirichlet algeba, the uniform algebra on the closed which is the uniform closure of the analytic polynomials on the unit circle. On the other hand the disk algebra on the closed unit disk,, that is, the uniform algebra of all complex-valued continuous functions on $\bar{D}$ which is analytic on the open unit disk $D$ is not a Dirichlet algebra. In fact, a uniform algebra needs not be Dirichlet in many cases. The ball algebra and the polydisk algebra on the ball and the polydisk of dimension 2 or greater are not Dirichlet algebra even on the \v Silov boundaries respectively. For further information see \cite{br,gamegame,stout}

According to \cite{dgpw} the numerical index of a Banach space was introduced by Lumer in 1968. For the algebra of all  bounded linear operator $L(B)$ on a Banach space $B$, the numerical index is 
\[
n(B)=\inf\{\nu(\bA):\bA\in L(B),\|\bA\|=1\},
\]
where $\nu(\bA)$ is the numerical radius
\[
\nu(\bA)=\sup\{|p(\bA(a))|: a\in S(B), p\in S(B^*), p(a)=1\}.
\]
Boyko, Kadets, Mart\'in and Mer\'i \cite[Theorem 2.1]{bkmm}
showed that the numetrical index of a lush space is $1$. Hence we see that 
\begin{crl}
Let $A$ be a uniform algebra. Then $n(A)=n(\overline{\operatorname{Re}A})=1$.
\end{crl}

\subsection{Extremely regular spaces and extremely $C$-regular spaces.}

An extremely regular space was given by Cengiz \cite{cengiz}.
An extremely $C$-regular space was introduced by Fleming and Jamison \cite[Definition 2.3.9]{fj1}.
\begin{dfn}
A $\K$-linear subspace $E$ of $C_0(Y,\K)$ is said to be {\em extremely $C$-regular (resp. regular)}  if for each $x$ in the Choquet boundary\footnote{The definition is given in Definition \ref{choquet}}  $\ch(E)$ (resp. $x\in Y$) satisfies the condition that for each $\varepsilon>0$ and each open neighborhood $U$ of $x$, there exists $f\in E$ such that $f(x)=1=\|f\|_\infty$, and $|f|<\varepsilon$ on $Y\setminus U$. 
\end{dfn}
Suppose that $m$ is a complex regular Borel continuous measure on $Y$. Then $E=\{f\in C_0(Y, \C):\int fdm=0\}$ is an extremely regular closed subspace of $C_0(Y,\C)$ (see \cite[Theorem]{cengiz}).
Another example is a closed subalgebra of $C_0(Y,\C)$ as follows.
Abrahamsen, Nygaard an P\~oldvere \cite{anp} introduced a {\em somewhat regular subspaces} of $C_0(Y,\K)$, which is a generalization of extremely regular subspaces.
\begin{dfn}[Definition 2.1 in \cite{anp}]
We call a $\K$-linear subspace $E$ of $C_0(Y,\K)$ {\em somewhat regular}, if for every non-empty open subset  $V$of $Y$ and $0<\varepsilon$, there exists $f\in E$ such that there exists $x_0\in V$ with $f(x_0)=1=\|f\|_\infty$ and $|f|\le \varepsilon$ on $Y\setminus V$.
\end{dfn}
\begin{prp}\label{alcreg}
A closed subalgebra of $C_0(Y,\C)$  which separates the points of $Y$ and has no common zeros is an extremely $C$-regular subspace of $C_0(Y,\C)$. 
In particular, a uniform algebra on a compact Hausdorff space $X$ is an extremely $C$-regular subspace of $C(X,\C)$.  If the Choquet boundary $\ch(A)$ is closed in $Y$ (resp. $X$), then $A|\ch(A)$ is an extremely regular subspace of $C_0(\ch(A),\C)$ if $\ch(A)$ is not compact and $C(\ch(A),\C)$ if $\ch(A)$ is compact. Let $S$ be the \v Silov boundary for $A$.
Then $A|S$ is a somewhat regular subspace of $C_0(S,\C)$ if $S$ is not compact and $C(S,\C)$ if $S$ is compact.
\end{prp}
\begin{proof}
The proof will be given after Corollary \ref{431}, where every thing will be ready.
\end{proof}

Abrahamsen, Nygaard and P\~oldvere \cite{anp} showed that 
extremely regular spaces play a role in recent theory of Banach spaces by
exhibiting that they involve the Daugave property, the symmetric strong diameter 2 property and so on 
 under some additional assumptions. 
We say $B$ has the {\em symmetric strong diameter 2 property} (SSD2P) if for every $\varepsilon>0$ and every finite collection of slices $S_1,\dots, S_m$, there 
exists $x_i\in S_i$ for $i=1,\dots, m$ and $y\in \ball(B)$ with $\|y\|>1-\varepsilon$ and $x_i\pm y\in S_i$ for all $i=1,\dots, m$.
Recall that a Banach space is {\em almost square} (ASQ) if whenever $x_1,\dots, x_n \in S(B)$, there exists a sequence $\{y_k\}\subset \ball(B)$ such that $\|x_i\pm y_k\|\to 1$ as $k\to \infty$ for all $1\le j\le n$.
Recall that a Banach space $B$ has the {\em Daugave property} if every rank-one operator $\bA$ on $B$ satisfies that $\|1+\bA\|=1+\|\bA\|$. 
Recall that a linear surjection $\mathcal{T}:N_1\to N_2$ for normed linear space $N_1$ and $N_2$ is called an $\varepsilon$-isometry if
\[
(1-\varepsilon)\|a\|\le \|\mathcal{T}(a)\|\le (1+\varepsilon)\|a\|
\]
for every $a\in N_1$.
We have
\begin{crl}\label{daugavet}
Let $A$ be a closed subalgebra of $C_0(Y,\C)$ which separates the points of $Y$ and has no common zeros.
Then we have
\begin{itemize}
\item[{\rm (1)}]
$A$ has the SSD2P,
\item[{\rm (2)}]
$A$ is ASQ if the \v Silov boundary is non-compact, 
\item[{\rm (3)}]
$A$ has the Daugavet property if the \v Silov boundary of $A$ is perfect,
\item[{\rm (4)}]
$A$ contains an $\varepsilon$-isometric copy of $c_0$ if $0<\varepsilon<1$.
\end{itemize}
\end{crl}
\begin{proof}
Put $A_0=A|S$, where $S$ is the \v Silov boundary. Then the restriction map is a surjective isometry from $A$ onto $A_0$.
The SSD2P, ASQ, the Daugavet propperty and to contain an $\varepsilon$-isometric copy of $c_0$ 
are inherited by an isometry, so it is enough to prove results for $A_0$.

By Proposition \ref{alcreg} that $A_0$ is a somewhat regular. Then by \cite[Theorems 2.2, 2.5, 2.6, 3.1]{anp} we have the conclusions.
\end{proof}
Note that Wojtaszczyk \cite[Theorem 2]{wo} proved that the Daugavet equation (DE) : $\|1+\bA\|=1+\|\bA\|$ holds for a weakly compact operator $\bA$ on a uniform algebra on  $X$ such that the strong boundary points are dense in $X$ and  $X$ has no isolated points.  As the strong boundary point coincides with the Choquet boundary points (Theorem \ref{chalg}) and they are dense in the \v  Silov boundary (Proposition \ref{431}), the hypothesis on the uniform algebra in the theorem of Wojtaszczyk  can be seen  that $A$ is a uniform algebra on a perfect $X$, where $X$ is the \v Silov boundary for $A$. As the $A$ and $A|S$ is isometric, where $A$ is a uniform algebra and $S$ is the \v Silov boundary for $A$, Wojtaszczyk in fact proved that the Daugavet equation holds for a weakly compact operators on a uniform algebra of which  \v Silov boundary is perfect.

At the end of the section we note 
that $C$-richness implies the somewhat regularity.
\begin{prp}
Let $E$ be a $\K$-linear subspace of $C_0(Y,\K)$. 
Suppose that for every nonempty open subset $U$ of $Y$ and $\varepsilon>0$, there exists a function $h_\varepsilon\in C_0(Y,\R)$ such that $0\le h_\varepsilon\le 1=\|h\|_\infty$  with support inside $U$ and that the $d(h_\varepsilon, E)<\varepsilon$.
Then $E$ is somewhat regular.
 We also have
\begin{itemize}
\item[{\rm (1)}]
$E$ has the SSD2P,
\item[{\rm (2)}]
$E$ is ASQ if $Y$ is non-compact, 
\item[{\rm (3)}]
$E$ has the Daugavet property if $Y$ is perfect,
\item[{\rm (4)}]
$E$ contains an $\varepsilon$-isometric copy of $c_0$ if $0<\varepsilon<1$.
\end{itemize}
In particular, if $Y$ is compact (hence $E$ is $C$-rich), then $E$ is somewhat regular and (1), (3) and (4) hold.
\end{prp}
\begin{proof}
Let $U$ be a nonempty open subset of $X$ and $\varepsilon>0$ arbitrary.
We prove that there exists $f\in E$ and $x_0\in U$ such that $f(x_0)=1=\|f\|$ and $|f|\le \varepsilon$. 
To prove it we may assume that $\varepsilon<1$. Put $\varepsilon_0=\varepsilon/(1+\varepsilon)$.
Then there exists $h\in C_0(Y,\R)$ such that $0\le h\le 1$ on $Y$, $h=0$ on $Y\setminus U$, $\|h\|_\infty=1$, and the distance between $h$ and $E$ is less than $\varepsilon_0$. Hence there exists $g\in E$ such that $\|h-g\|_\infty<\varepsilon_0$. Then $\|g\|_\infty>\|h\|_\infty-\varepsilon_0=1-\varepsilon_0$. As $h=0$ on $Y\setminus U$, we infer that $|g|<\varepsilon_0$ on $Y\setminus U$. Choose $x_0\in Y$ so that $\|g\|_\infty=|g(x_0)|$. 
Put $f=g/g(x_0)$. Then $f\in E$, $f(x_0)=1=\|f\|_\infty$ and $|f|<\varepsilon_0/(1-\varepsilon_0)=\varepsilon$ on $Y\setminus U$. As $U$ and $\varepsilon$ are arbitrary, we have that $E$ is somewhat regular.
Then by \cite[Theorems 2.2, 2.6, 3.1]{anp}  we have the coclusion.
\end{proof}
 

\section{Strong boundary points, Choquet boundary points and \v Silov boundary points}\label{sec3}


\subsection{Function spaces which  strongly separate the points in the underlying spaces}\label{sec X}
The definition of ``strongly separate the points''  in this paper is due to that of Araujo and Font \cite{af97}. Some authors use this term  in a different notion (cf. Stout \cite[p.36]{stout} and Miura \cite[p.779]{miuracent} )
\begin{dfn}
Let $E$ be a complex (resp. real) linear subspace of $C_0(Y,\C)$ (resp. $C_0(Y,\R)$). 
We say that $E$ (resp. strongly)  separates the points of $Y$, if for every pair $x,y\in Y$ with $x\ne y$ there exists $f\in E$ such that $f(x)\ne f(y)$ (resp. $|f(x)|\ne |f(y)|$). We say that $E$ has no common zeros if for every $x\in Y$ there exists $f\in E$ such that $f(x)\ne 0$. 
\end{dfn}
The following example exhibits an important space of functions which separates, but does not strongly separate the points of the underlying space. 
\begin{exa}
The space $C_0(Y,\R)$ and $C_0(Y,\C)$ strongly separate points of $Y$ and have no common zeros by the Urysohn's lemma. Let $B$ be a complex Banach space. For $a\in B$ we denotes $\hat a: \ball(B^*)\setminus \{0\}\to \C$ by $\hat{a}(q)=q(a)$ for $q\in \ball(B^*)\setminus \{0\}$. Then $\widehat{B}=\{\hat{a}: a\in B\}$ is a uniformly closed subspace of $C_0(\ball(B^*)\setminus \{0\},\C)$ which separates the points of $\ball(B^*)\setminus \{0\}$. On the other hand $|\hat{a}(p)|=|\hat{a}(-p)|$ for any $a\in B$ and $p\in \ball(B^*)\setminus \{0\}$, that is, $\widehat{B}$ does not {\em strongly} separate the points of the underlying space $\ball(B^*)\setminus\{0\}$. The situation is similar for a real Banach space.
\end{exa}
If $E$ contains constants (in this case it is automatic that $Y$ is compact), then $E$ separates the points of $Y$ if and only if $E$ strongly separates the  points of $Y$.
A subalgebra of $C_0(Y,\C)$ separates the points of $Y$ if and  only if $A$ strongly separates the points. In fact, we have the following.
\begin{prp}\label{ssp}
Supposes that $E$ is a $\K$-linear subspace of $C(X,\K)$,  such that $1\in E$, where $X$ is a compact Hausdorff space. If $E$ separates the points of $X$, then $E$ strongly separates the points of $X$.
Suppose that $A$ is a  subalgebra of $C_0(Y,\K)$ which separates the points of $Y$. Then $A$ strongly separates the points of $Y$.
\end{prp}
\begin{proof}
We prove the first assertion. Suppose that $x,y\in X$ with $x\ne y$. Then there exists $f\in E$ such that $f(x)\ne f(y)$. If $f(x)=0$, then $|f(x)|\ne |f(y)|$. If $f(x)\ne 0$ then $|f|$ or $|f+f(x)|$ separates $x$ and $y$, where $f+f(x)\in E$ since $E$ contains constants.

Next suppose that  $A$ is a subalgebra of $C_0(Y,\K)$.
Suppose that $x,y\in Y$ with $x\ne y$. Then there exists $f\in A$ such that $f(x)\ne f(y)$. If $f(y)=0$, then $f(x)\ne 0$ and $|f(x)|\ne |f(y)|$ follows. We may assume that $f(y)=1$.  Suppose that $|f(x)|=|f(y)|$. Then there exists a complex number $\lambda$ with unit modulus such that $f(y)=\lambda f(x)$. As $\lambda\ne 1$, we infer that $|\lambda^2+\lambda|<2$. Hence  we have $|f(x)+f(x)^2|=|\lambda+\lambda^2|<1$ and 
$2=|f(y)+f(y)^2|$. Hence $f+f^2\in A$ strongly separates $x$ and $y$.
\end{proof}


\subsection{Strong boundary points}\label{sec X}
We recall the notion of a {\em strong boundary point} for a certain algebra of functions from \cite[Definition 7.6]{stout}. We adapt it to  subspaces of $C_0(Y,\K)$. 
Note that the two definitions of ours (Definition \ref{4.3}) and that  in \cite[Definition 2.3.9]{fj1}  are equivalent if the space is closed under the multiplication while it is not in general.
\begin{dfn}\label{4.3}
Let $E$ be a $\K$-linear subspace of $C_0(Y,\K)$.
We say that a point $x_0\in Y$ is a {\em strong boundary point} of $E$ if for each open neighborhood $U$ of $x_0$ there exists $f\in E$ such that $1=f(x_0)=\|f\|_\infty$, and $|f(x)|<1$ for all $x\in Y\setminus U$. 
We say that a closed subset $K$ of $Y$ is a {\em peak set} for $E$ if there exists a function $f\in E$ such that $f=1$ on $K$ and $|f|<1$ on $Y\setminus K$. We say that any such  $f$ {\em peaks} on $K$ and $f$ is a {\em peaking function} for $K$. A {\em weak peak set} (or {\em peak set in the weak sense}) for $E$ is a non-empty intersection of peak sets for $E$. A point $y_0\in Y$ is called a {peak point} for $E$ if $\{y_0\}$ is a peak set for $E$. A point $y_0\in Y$ is called a {\em weak peak point} (or {\em peak point in the weak sense}) for $E$ if $\{y_0\}$ is a weak peak set  for $E$.
\end{dfn}

Recall that the peripheral range $\rpi(f)$ of $f\in C_0(Y,\C)$ is the set $\{z\in f(Y): |z|=\|f\|_\infty\}$.
A function $f\in E$ is a peaking function for a closed subsest of $Y$ if and only if $\rpi(f)=\{1\}$.
\begin{prp}\label{4.4}
Let $E$ be a uniformly closed $\K$-linear subspace of $C_0(Y,\K)$. Let $x_0\in Y$.
Then the following are equivalent.
\begin{itemize}
\item[\rm{(i)}]
The point $x_0$ is a strong boundary point for $E$,
\item[\rm{(ii)}]
for every open neighborhood $U$ of $x_0$, there exists $f\in E$ such that $f(x_0)=1=\|f\|_\infty$ and $|f(x)|<1$ for every $x\in Y\setminus U$ and  $\rpi(f)=\{1\}$,
\item[\rm{(iii)}]
the point $x_0$ is a weak peak point for $E$
\end{itemize}
\end{prp}
\begin{proof}
We prove (i) implies (ii). Suppose (i) holds. Let $U$ be an open neighborhood of $x_0$. Then there exists $f\in E$ with $f(x_0)=1=\|f\|_\infty$ and $|f|<1$ on $Y\setminus U$.  Put $U_n=\{y\in Y:|f(y)-1|<1/n\}$ for each positive integer $n$. Then $U\cap U_n$ is an open neighborhood of $x_0$. As $f^{-1}(1)=\bigcap_{n=1}^\infty U_n$, $f^{-1}(1)$ is a $G_{\delta}$ set. By (i) there exists a function $f_n\in E$ such that $f_n(x_0)=1=\|f_n\|_\infty$, $|f_n|<1$ on $Y\setminus (U\cap U_n)$. Put $g=\sum_{n=1}^\infty(f_n/2^n)$. Then $g\in E$ since $E$ is closed in $C_0(Y, \K)$. For every $x\in Y\setminus  U$, we have $x\in Y\setminus (U\cap U_1)$, hence $|f_1(x)|<1$, so $|g(x)|<1$. Suppose that $|g(y)|=1$ for some $y\in Y$. Since $1=|g(y)|\le \sum(|f_n(y)|/2^n\le 1$, we infer that $|f_n(y)|=1$ for every $n$. We conclude that $x\in \bigcap_{n=1}^\infty U_n$. Put $h=(f+g)/2$.  Then $h(x_0)=1=\|h\|_\infty$ and $|h|<1$ on $Y\setminus U$. Suppose that $|h(z)|=1$ for  $z\in Y$. Since 
\begin{equation}\label{hfg}
1=|h(z)|=|(f(z)+g(z)|/2\le (|f(z)|+|g(z)|)/2\le 1
\end{equation}
and $\|f\|_\infty=\|g\|_\infty=1$,  we infer that $|g(z)|=1$. Hence $z\in \bigcap_{n=1}^\infty U_n$. Since $f^{-1}(1)=\bigcap_{n=1}^\infty U_n$ we infer that $f(z)=1$. By \eqref{hfg} we get $g(z)=1$. Therefore $h(z)=1$. Thus $\rpi(h)=\{1\}$.

To prove (ii) implies (iii), assume  (ii). Let $y\in Y\setminus\{x_0\}$. There is an open neighborhood $U_y$ of $x_0$ such that $y\not\in U_y$. By (ii) there is a function $f_y\in E$ such that $f_y(x_0)=1=\|f_y\|$, $|f_y|<1$ on the closed set $Y\setminus U_y$ , and $\rpi(f_y)=\{1\}$. Then $f_y^{-1}(1)$ is a peak set which contains $x_0$ and $y\not\in f_y^{-1}(1)$. Hence $\bigcap_{y\in Y\setminus \{x_0\}}f_y^{-1}(1)=\{x_0\}$ is a weak peak set, so $x_0$ is a weak peak point for $E$.

To prove (iii) implies (i), assume  (iii): there exists a family of peak sets $\{K_{\alpha}\}$ such that $\bigcap_{\alpha}K_{\alpha}=\{x_0\}$. Suppose that $U$ is an open neighborhood of $\{x_0\}$. Then $\bigcap_{\alpha}K_{\alpha}\subset U$. By considering the one point compactification, if necessary, we infer that there is a finite number of $\alpha_1,\dots, \alpha_n$ such that $\bigcap_{k=1}^nK_{\alpha_k}\subset U$. Let $f_j\in E$ be a function which peaks on $K_{\alpha_j}$. Then $f=\frac{1}{n}\sum_{j=1}^nf_j$ is in $E$ and peaks on $\bigcap_{j=1}^nK_{\alpha_j}$. Hence $f(x_0)=1=\|f\|_\infty$. Since $Y\setminus U\subset Y\setminus \bigcap_{j=1}^nK_{\alpha_j}$, we have $|f|<1$ on $Y\setminus U$.
\end{proof}
Note that if $E$ is a $\K$-linear subspace of $C(X,\K)$ for a compact Hausdorff space $X$ and $1\in E$, which needs not to be uniformly closed, then (i) of Proposition \ref{4.4} implies (ii) of Proposition \ref{4.4}. In fact, if $x_0\in X$ is a strong boundary point and $U$ is an open neighborhood and $f\in E$ satisfies $1=f(x_0)=\|f\|_\infty$. Then $(f+1)/2$ satisfies the condition (ii).
On the other hand if $E$ is not uniformly closed nor $1\not\in E$, then (i) needs not imply (ii).
\begin{exa}
Let $Y=\{z\in \C:|z|=1\}$. For any positive integer $n$, let $f_n$ be a continuous function on $Y$ such that 
$f_n(z)=z$ for any $z\in Y$ with $|1-z|\le 1/n$ and $|f_n(z)|<1$ for any $z\in Y$ with $|1-z|>1/n$. Let $E$ be a linear subspace generated by $\{1\}\cup\{f_n\}_{n=1}^\infty$.  Note that $E$ is a $\K$-linear subspace of $C(Y,\K)$ which is not uniformly closed. Then $1$ is a strong boundary point for $E$, while $1$ does not satisfy the condition (ii) of Proposition \ref{4.4}
\end{exa}
Note that the definition of a strong boundary point sometimes differs from authors to authors. Stout in \cite[Definition 7.6]{stout} sais that $x_0\in Y$ is a strong boundary point for a certain subalgebra of $C_0(Y,\C)$ with some condition if for each open neighborhood $U$ of $x_0$ there exists $f\in A$ such that $1=f(x_0)=\|f\|_\infty$ and $|f(x)|<1$ for all $x\in Y\setminus U$.  Aroujo and Font \cite[p. 80]{af}  follow this definition not only for a subalgebra but also for a linear subspace.     
Rao and Roy \cite[Definition 8]{rr} recall the notion of a strong boundary point from \cite{stout} for a linear subspace. 
Note that a weak peak point (peak point in the weak sense in \cite{br})  for a uniform algebra (function algebra in \cite{br}) is also referred to as a strong boundary point in \cite[p.96]{br}. 
In a book of Gamelin \cite{gamegame} a weak peak set (resp. point) is referred to as a $p$-set (resp. point) or  a generalized peak set (resp. point).

The definition of a strong boundary point for a $\K$-linear subspace $E$ of $C_0(Y,\K)$ due to Fleming and Jamison \cite[Definition 2.3.9]{fj1}: a point $x_0\in Y$ is a strong boundary point of $E$ if for each open neighborhood $U$ of $x_0$, and each $\varepsilon>0$, there exists $f\in E$ such that $1=f(x_0)=\|f\|$, and $|f(x)|<\varepsilon$ for all $x\in Y\setminus U$, is stronger than that given in Definition \ref{4.3}. 

If the corresponding space $E$ is an algebra, it is easy to see that two definitions (one by Fleming and Jamison \cite{fj1}, and one by Definition \ref{4.3}) on the strong boundary point coincide with each other. As a corollary of Proposition \ref{4.4} we have the following.

\begin{crl}\label{37}
Let $A$ be a closed subalgebra of $C_0(Y,\K)$. Let $x_0\in Y$. Then the following are equivalent.
\begin{itemize}
\item[\rm{(i)}] 
The point $x_0$ is a strong boundary point for $A$,
\item[\rm{(ii)}]
For every open neighborhood $U$ of $x_0$, and $\varepsilon>0$, there exists $f\in A$ such that $f(x_0)=1=\|f\|_\infty$ and $|f(x)|<\varepsilon$ for every $x\in Y\setminus U$,
\item[\rm{(iii)}]
For every open neighborhood $U$ of $x_0$, and $\varepsilon>0$, there exists $f\in A$ such that $f(x_0)=1=\|f\|_\infty$ and $|f(x)|<\varepsilon$ for every $x\in Y\setminus U$ and  $\rpi(f)=\{1\}$,
\item[\rm{(iv)}]
The point $x_0$ is a weak peak point for $A$
\end{itemize}
\end{crl}
\begin{proof}
The point of the proof is as follows. Suppose that  $x_0\in Y$. If $x_0\in U$ is open and $f\in A$ satisfies  $1=f(x_0)=\|f_0\|_\infty$ and $|f|<1$ on $Y\setminus U$, then for every $\varepsilon>0$ there exists a positive interger such that $|f^n|<\varepsilon$ on $Y\setminus U$. It is true since $Y\setminus U$ is closed and hence there is $r<1$ such that $|f|<r$ on $Y\setminus U$. 
After that, a routine argument proves the corollary. We omit a detailed proof.
\end{proof}
Corollary \ref{37} does not hold if $E$ is not closed under multiplication.  
\begin{exa}\label{0ofE}
Let $f_0: (0,1]\to \R$ be defined as 
\begin{equation*}
f_0(t)=
\begin{cases}
t,& 0<t\le\frac12, 
\\
-\frac32 t+\frac54, & \frac 12\le t\le 1.
\end{cases}
\end{equation*}
Let $g:(0,1]\to \R$ be
\begin{equation*}
g_0(t)=
\begin{cases}
0, &0<t\le\frac12 
\\
-(t-\frac12)(t-1), &\frac12\le t\le 1.
\end{cases}
\end{equation*}
Put $E=\{\lambda f_0+\mu g_0:\lambda, \mu \in \K\}$. Then $E$ is a $\K$-linear subspace of $C_0((0,1],\K)$ which strongly separates the points of $Y$. 
The point $\frac12$ is a strong boundary point (in the sense of Definition \ref{4.3}) for the space $E$ 
while it does not satisfy the condition due to Fleming and Jamison \cite[Definition 2.3.9]{fj1}.
\end{exa}
To study subspace  $E$ in $C_0(Y,\K)$ it is usufull to consider the addition of constant.
We add constant functions in a subspace of $C_0(Y,\K)$.
\begin{dfn}
For a locally compact Hausdorff space $Y$, we denote by $Y_{\infty}=Y\cup\{\infty\}$ the one-point-compactification of $Y$. Let $E$ be a $\K$-linear subspace of $C_0(Y,\K)$. For $f\in E$ we denote by $\dot{f}$ the unique extension of $f$ on $Y$,
\begin{equation*}
\dot{f}(y)=
\begin{cases}
f(y), & y\in Y \\
0, & y=\infty.
\end{cases}
\end{equation*}
Then $\dot{f}$ is continuous on $Y_\infty$. We denote
\[
\dot{E}+\K=\{F\in C(Y_\infty,\K): F=\dot{f}+c, f\in E, c\in \K\}.
\]
Then  $\dot{E}+\K$ is a $\K$-linear subspace of $C(Y_\infty,\K)$ which contains constants. 
\end{dfn}
We may sometimes suppose that $E$ is a closed subspace of $\dot{E}+\K$ without a confusion.
It is easy to see that  $\dot{E}+\K$ separates the points of $Y_\infty$ and has no common zeros provided that  $E$ separates the points of $Y$  and it has no common zeros. 
By a routine work we have that  $\dot{E}+\K$ is closed in $C(Y_\infty,\K)$
if $E$ is closed in $C_0(Y,\K)$.
It is also a routine work to see that for $F\in \dot{E}+\K$, $F=\dot{f}$ for a $f\in E$ if and only if $F(\infty)=0$.
\begin{lmm}\label{4.7}
Suppose that $E$ is a $\K$-linear subspace of $C_0(Y,\K)$. If $x_0\in Y$ is a strong boundary point for $E$, then $x_0$ is a strong boundary point for $\dot{E}+\K$.
\end{lmm}
\begin{proof}
The proof is trivial and is omitted.
\end{proof}
The converse of the above lemma does not hold in general. 
\begin{exa}\label{1ofE}
Let $E$ be the same space defined in Example \ref{0ofE}.
Then $1$ is not a strong boundary point for $E$. On the other hand, $1$ is a strong boundary point for $\dot{E}+\K$.
\end{exa}
The converse of Lemma \ref{4.7} holds for a closed subalgebra of $C_0(Y,\K)$.
\begin{prp}\label{49}
Suppose that $A$ is a closed subalgebra of $C_0(Y,\K)$. A point $x_0\in Y$ is a strong boundary point for $A$ if and only if it is a strong boundary point for $\dot{A}+\K$.
\end{prp}
\begin{proof}
Suppose that $x_0\in Y$ is a strong boundary point for $\dot{A}+\K$. For any open neighborhood $U$ of $x_0$ in $Y$, $U$ may be considered as an open neighborhood of $x_0$ in $Y_\infty$. Hence there exists a function $F\in \dot{A}+\K$ such that $F(x_0)=1=\|F\|_\infty$ and $|F|<1$ on $Y_\infty\setminus U$. We note that $\infty \not\in U$. Hence we may suppose that $|F(\infty)|<1$. Put 
$\pi :\{z\in \C:|z|\le 1\} \to \{z\in \C:|z|\le 1\}$ by 
\[
\pi(z)=\frac{1-\overline{F(\infty)}}{1-F(\infty)}\cdot\frac{z-F(\infty)}{1-\overline{F(\infty)}z},\quad z\in \{z\in \C:|z|\le 1\}
\]
if $\K=\C$. We infer that $\pi(F(\infty))=0$ and $\pi(1)=1$.
As $\pi$ is uniformly approximated by analytic polynomials (in fact, $\pi(rz)$ for any $0<r<1$ is uniformly  approximated by the Taylor expansion on $\{z\in \C:|z|\le 1\}$, and $\pi(rz)\to \pi(z)$ uniformly on $\{z\in \C:|z|\le 1\}$) and as $\dot{A}+\C$ is uniformly closed algebra,  we have $\pi\circ F\in \dot{A}+\C$.
If $\K=\R$, then put
 $\pi: [-1,1]\to [-1,1]$ by 
\begin{equation*}
\pi(t)=
\begin{cases}
\frac{1}{1+F(\infty)}t-\frac{F(\infty)}{1+F(\infty)}, & -1\le t\le F(\infty)
\\
\frac{1}{1-F(\infty)}t-\frac{F(\infty)}{1-F(\infty)}, &F(\infty)\le t \le 1.
\end{cases}
\end{equation*}
Then $\pi(F(\infty))=0$ and $\pi(1)=1$. 
By the Weierstrass approximation theorem $\pi$ is uniformly approximated by polynomials on $[-1,1]$. Hence $\pi\circ F\in \dot{A}+\R$. 
In any case we have that $\pi\circ  F\in \dot{A}+\K$, $\pi\circ F(\infty)=0$. 
Thus there exists $f\in A$ such that $F=\dot{f}$. Then we have that $\|f\|_\infty=\|F\|_\infty=1$ and $f(x)=\pi\circ F(x)=\pi(1)=1$. As $|\pi(z)|<1$ for any $|z|<1$, we have $|\pi\circ F|<1$ on $Y_\infty \setminus U$. Therefore $|f|<1$ on $Y\setminus U$. It follows that $x_0$ is a strong boundary point for $A$.

The converse statement is just Lemma \ref{4.7}
\end{proof}


\subsection{The Choquet boundary and the \v Silov boundary}
The Choquet boundary was first mentioned by that name in a paper of Bishop and de Leeuw \cite[p.306]{bd}.
Definitions of Choquet boundary are differs from case by case, although they are equivalent in the possible situation where the difinitions can be applied.
A definition of the Choquet boundary for $\K$-linear subspace of $C(X,\K)$, for a compact Hausdorff space $X$, which contains constant functions is described by Phelps in \cite[Section 6]{phelps}. 
Let $E$ be a $\K$-linear subspace of $C(X,\K)$. Suppose that $1\in E$. The state space of $E$ is 
$K(E)=\{\phi\in E^*: \phi(1)=1=\|\phi\|\}$ (cf. \cite[p.27]{phelps}).
The Choquet boundary in \cite{phelps} is defined as the set of all $x\in Y$ such that  the point evaluation $\tau_x$ is in $\ext(K(E))$, the set of all extreme points of the state space. Since it is easy to see that $\T\ext(K(E))=\ext(B(E^*))$, hence the two definitions of the Choquet boundary (Definition \ref{choquet} and Definition in \cite[p.29]{phelps})  are equivalent provided that $Y$ is compact and $1\in E$.
In section 8 of \cite{phelps} Phelps describes an equivalent form of the Choquet boundary for a uniform algebra in terms of strong boundary points by referring a theorem of Bishop and de Leeuw. Browder \cite[Section 2-2]{br} exhibits a definition of the Choquet boundary for a uniform algebra in terms of mesures, which is also equivalent to that exhibited in  Definition \ref{choquet} in the case of a uniform algebra.  Rao and Roy \cite[p.176]{rr} defines the Choquet boundary for a uniformly closed complex linear subspace of $C(X,\C)$ which separates the points of a compact Hausdorff space $X$, in a similar way as our Definition \ref{choquet}.

We expect that some  results in this subsection are a part of folklore, but for the sake of a self-contained exposition we have included as possible as complete proofs of all the results stated.

\subsubsection{Definition of the Choquet boundary.}

We recall the notion of the {\em Choquet boundary} for a $\K$-subspace of $C_0(Y,\K)$ from \cite[Definition 2.3.7]{fj1}, which was stated by Novinger \cite[p.274]{nov}. See also \cite{semad}.

For a $\K$-linear subspace $E$ of $C_0(Y,\K)$ and $x\in Y$, $\tau_x$ denotes the point evaluation at $x$, that is,  $\tau_x:E\to \K$ such that $\tau_x(f)=f(x)$ for $f\in E$.
If $E$ contains constant, then $\|\tau_x\|=1$. 
In general,  $\tau_x\in \ball(E^*)$ and $\|\tau_x\|$ needs not be $1$. For example, put $E=\{f\in P(\bar{D}):f(0)=0\}|(\bar{D}\setminus\{0\})$, where $P(\bar{D})$ is the disk algebra on the closed unit disk $\bar{D}$ in the complex plane. Let $x\in D$, the open disk. Then $\|\tau_x\|=|x|$ by the Schwarz lemma.  

We define the Choquet boundary for a $\K$-linear subspace which needs not to be closed, not to separates the points of the underlying space, may have common zeros.
\begin{dfn}\label{choquet}
Suppose that $E$ is a $\K$-linear subspace of $C_0(Y,\K)$. The Choquet boundary for $E$ denoted by $\ch(E)$ is the set of all $x\in Y$ such that the point evaluation $\tau_x$ is in $\ext(\ball(E^*))$, the set of all extreme points of $\ball (E^*)$.
\end{dfn}


\subsubsection{
The representing measures and 
the Arens-Kelley theorem revisited}

It is crucial for the foregoing discussion on the Choquet boundaries that measure theoretic arguments should be concerned.
We begin by recalling the representing measures for bounded linear functionals.
\begin{dfn}\label{representingm}
Let $E$ be a $\K$-linear subspace of $C_0(Y,\K)$. Suppose that $\phi\in E^*$. We say that a complex regular Borel measure $m$ on $Y$ is a {\em representing measure} for $\phi$ if
\[
\phi(f)=\int fdm,\quad f\in E
\]
and $\|m\|=\|\phi\|$, where $\|m\|=|m|(Y)$ is the total valuation of $m$.
\end{dfn}
Existence of a representing measure for any $\phi\in E^*$ is as follows: let $\phi\in E^*$. By the Hahn-Banach extension theorem there is $\Phi\in C_0(Y,\K)^*$ which extends $\phi$ and $\|\phi\|=\|\Phi\|$. Then by the Riesz-Kakutani theorem, there exists a complex regular Borel measure $m$ such that $\Phi(f)=\int fdm$ for every $f\in C_0(Y,\K)$ and $\|m\|=\|\Phi\|$. Hence $m$ is a representing measure for $\phi$. Recall that the {\em support} $\suppm(m)$of a complex regular Borel measure $m$ is the set
\[
\{x\in Y: \text{$|m|(G)>0$ for every open neighborhood $G$ of $x$}\},
\]
where $|m|$ is the total valuation measure of $m$.

Versions of the Arens-Kelley theorem, which characterize extreme points of the unit ball in the dual space of a subspace of $C_(Y,\K)$,  have been obtained by a variety of authors. 
The Arens-Kelley theorem and the following corollary are well known, however for the sake of a self-contained exposition and for the convenience of the readers we include complete proofs.

For $x\in Y$ we denote the {\em point mass} at $x$ by $D_x$: $D_x$ is a complex regular Borel measure on $Y$ such that $D_x(\{x\})=1=\|D_x\|$. 
\begin{thmak}
Suppose that $\phi\in \ext\ball(C_0(Y,\K)^*)$. Then there exists a unique $x\in Y$ and $\lambda\in \T$ such that $\phi=\lambda \tau_x$. The representing measure for $\phi$ is only $\lambda D_x$.
Conversely, $\lambda \tau_x$ is an extreme point of $\ball(C_0(Y,\K)^*)$ for every $x\in Y$ and $\lambda\in \T$.
\end{thmak}
\begin{proof}
Let $\phi\in \ext\ball(C_0(Y,\K)^*)$ and $m$ a representing measure for $\phi$. 
Let $y\in \suppm(m)$ arbitrary. Suppose that $|m|(U)=1$ for every open neighborhood $U$  of $y$. 
By the regularity of $m$ we infer that $|m|(\{y\})=1=\|m\|$, hence $m=\lambda D_y$ for unimodular complex number $\lambda$. 

Suppose that there exists an open neighborhood $U_0$ of $y$ with $|m|(U_0)<1$. As $y\in \suppm(m)$, $0<|m|(U_0)$ holds. 
By the definition of $|m|$ we have 
\[
|m|(U_0)=\sup\{\sum |m(G_j)|: \text{$G_j$ is a Borel set, $\cup_j G_j=U_0$, $G_i\cap G_j=\emptyset$ for $i\ne j$}\},
\] hence there is $G_j$ such that $m(G_j)\ne 0$. 
By the (inner) regularity of $m$, there exists a compact subset $L\subset G_j\subset U_0$ with $0<|m(L)|$.  Also there exists an open set $V$ with $L\subset V\subset U_0$ such that $|m|(V\setminus L)<|m(L)|/2$. By the Urysohn's lemma there exists $f_0\in C_0(Y,\R)$ such that $0\le f_0\le 1$, $f_0=1$ on $L$, and $f_0=0$ on $Y\setminus V$.
Then we have 
\[
0<|m(L)|\le |m|(L)\le |m|(V)\le |m|(U_0)<1.
\]
Hence we also have $0<|m|(Y\setminus V)<1$.
Put
\[
\phi_1(g)=\frac{1}{|m|(V)}\int_V gdm,\quad g\in C_0(Y,\K)
\]
and 
\[
\phi_2(g)=\frac{1}{|m|(Y\setminus V)}\int_{Y\setminus V}g dm,\quad g\in C_0(Y, \K).
\]
As $m$ is a representing measure for $\phi$ we have  $\phi=|m|(V)\phi_1+|m|(Y\setminus V)\phi_2$, where  $|m|(V)+|m|(Y\setminus V)=\|m\|=1$. As $\phi\in \ext\ball (C_0(Y, \K)^*)$ we have $\phi=\phi_1=\phi_2$. Then 
\[
\phi (f_0)=\phi_2(f_0)=\frac{1}{|m|(Y\setminus V)}\int_{Y\setminus V}f_0 dm=0
\]
since $f_0=0$ on $Y\setminus V$. On the other hand
\begin{multline*}
|\phi(f_0)|=|\phi_1(f_0)|=\left|\frac{1}{|m|(V)}\int_V f_0 dm\right|
\ge
\frac{1}{|m|(V)}\left(|\int_Lf_0 dm|-|\int_{V\setminus L}f_0 dm|\right)
\\
\ge
\frac{1}{|m|(V)}\left(|m(L)|-|m|(V\setminus L)\right)
>
\frac{|m(L)|}{2|m|(V)}>0,
 \end{multline*}
hence $\phi(f_0)\ne 0$, which is a cotradiction.

Suppose that $m$ is a representing measure for $\lambda\tau_x$.
Let $U$ be an open neighborhood of $x$. Then by the Urysohn's lemma there exists $f\in C_0(Y,\K)$ such that $0\le f\le 1$ on $Y$, $f(x)=1$, and $f=0$ on $Y\setminus U$. Since 
\[
1=f(x)=\tau_x(f)=|\int f dm|\le \int_U|f|d|m|\le|m|(U)\le 1,
\]
we see that $\suppm(m)\subset U$. As $U$ can be arbitrary, we see that $\suppm(m)=\{x\}$. Thus we infer that $m=\lambda D_x$.

Suppose conversely that $x_0\in Y$ and $\lambda_0\in \T$. 
Suppose that $\phi_0=(\phi_1+\phi_2)/2$ for $\phi_j\in \ball(C_0(Y,\K)^*)$, $j=1,2$. Let $\mu_1$ be a representing measure for $\phi_1$. We prove that $\suppm(\mu_1)=\{x_0\}$. Suppose not. As the support of a regular measure is not empty, there exists $y\in \suppm(\mu_1)\setminus \{x_0\}$. Then by the Urysohn's lemma there exists $f\in C_0(Y,\R)$ such that $f(y)=0\le f\le 1=f(x_0)$. Put $U=\{z\in Y: f(z)<1/2\}$. Then $U$ is an open neighborhood of $y$,  and $0<|\mu_1|(U)<1$ since $y\in \suppm(\mu_1)$. Thus
\[
|\phi_1(f)|\le \int_U|f|d|\mu_1|+\int_{Y\setminus U}|f|d|\mu_1|\le\frac12|\mu_1|(U)+|\mu_1|(Y\setminus U)<1.
\]
It follows that 
\[
1=|\phi_0(f)|\le (|\phi_1(f)|+|\phi_2(f)|)/2 <1,
\]
which is a contradiction proving that $\suppm(\mu_1)= \{x_0\}$. We infer that   $\phi_1=\lambda_1\tau_{x_0}$ for $\lambda_1\in \T$. In the same way we have that $\phi_2=\lambda_2\tau_{x_0}$ for $\lambda_2\in \T$.
Since
\[
\lambda_0=\phi_0(f)=(\phi_1(f)+\phi_2(f))/2=(\lambda_1+\lambda_2)/2
\]
and $|\lambda_j|=1$ for $j=0,1,2$ we infer that $\lambda_0=\lambda_1=\lambda_2$.  Thus $\phi_1=\phi_2=\phi_0$. We concluded that $\phi_0\in \ext\ball(C_0(Y,\K)^*)$.
\end{proof}
\begin{crl}\label{akc}
Suppose that $E$ is a $\K$-linear subspace of $C_0(Y,\K)$. Suppose that $\phi\in \ext\ball(E^*)$. Then there exists $y\in Y$ and $\lambda\in \T$ such that $\phi=\lambda\tau_y$.
\end{crl}
\begin{proof}
Put
\[
S=\{\varphi\in C_0(Y,\K)^*: \text{$\varphi$ is a Hahn-Banach extension of $\phi$}\}.
\]
Then $S$ is a non-empty weak$^*$-closed convex subset of $C_0(Y,\K)$. Then the Krein-Milman theorem asserts that there exists a $\Phi\in \ext(S)$. Then $\Phi$ is an extreme point of $\ball(C_0(Y,\K)^*)$. In fact, suppose that $\Phi=(\Phi_1+\Phi_2)/2$ for $\Phi_1,\Phi_2 \in \ball(C_0(Y,\K)^*)$. Then  $\phi=\Phi|E=(\Phi_1|E+\Phi_2|E)/2$, and $\Phi_j|E\in \ball(E^*)$ for $j=1,2$. 
As $\phi$ is an extreme point of $\ball(E^*)$, we have  $\Phi_1|E=\Phi_2|E=\phi$. Since $1=\|\phi\|=\|\Phi_j|E\|\le \|\Phi_j\|=1$, we have $\|\Phi_j\|=1$ for $j=1,2$. 
Hence $\Phi_j\in S$ for $j=1,2$. As $\Phi$ is an extreme point in $S$, we have $\Phi=\Phi_1=\Phi_2$. Thus $\Phi\in \ext(C_0(Y,\K)^*)$. By the Arens-Kelley theorem there exists $y\in Y$ and $\lambda\in \T$ such that $\Phi=\lambda D_y$. Thus we see that $\phi=\Phi|E= \lambda D_y$.
\end{proof}
Note that $\lambda\tau_x$ needs not to be an extreme point of $\ball(E^*)$ in general. In fact, $\tau_0\not\in \ext\ball(P(\bar{D})^*)$ for the disk algebra on the closed unit disk $\bar{D}$ in the complex plane.

\subsubsection{
The Choquet boundary points, uniqueness of representing measures and strong boundary points.}

Let $E$ be a $\K$-linear subspace of $C_0(Y,\K)$ and $x\in Y$. We study the relationship of the following (i), (ii) and (iii) :
\begin{itemize}
\item[{\rm(i)}]
$x$ is a strong boundary point for $E$,
\item[{\rm(ii)}]
the representing measure for $\tau_x$ is only $D_x$,
\item[{\rm{(iii)}}]
$x\in \ch(E)$.
\end{itemize}

We recall the definition of the boundary.
\begin{dfn}
Suppose that $E$ is a $\K$-linear subspace of $C_0(Y,\K)$. A subset $L$ of $Y$ is said to be a boundary if for each $f\in E$ there exists a point $x\in L$ such that $|f(x)|=\|f\|_\infty$. 
\end{dfn}
The following may be well known, for example, \cite[Proposition 6.3]{phelps} states about the case where $E$ contains $1$ and $Y$ is compact. However for the sake of a self-contained exposition and for the convenience of the readers we included a complete proof.
\begin{prp}\label{chbou}
Suppose  that $E$ is a $\K$-linear subspace of $C_0(Y,\K)$. Then the Choquet boundary $\ch(E)$ is a boundary.
\end{prp}
\begin{proof}
Let $f\in E$. We may assume that $\|f\|_\infty=1$. 
Then there exists $y\in Y$ such that $|f(y)|=1$. Put $L=\{\phi\in \ball(E^*):\phi(f)=f(y)\}$. As $\tau_y\in L$, $L$ is non-empty weak*-closed convex subset of $\ball(E^*)$. The Krein-Milman theorem asserts that there exists $\phi_0\in \ext L$. Then $\phi_0$ is an extreme point of  $\ball(E^*)$. In fact, suppose that $\phi_0=(\phi_1+\phi_2)/2$ for $\phi_1,\phi_2\in \ball(E^*)$. Then by
\[
1=|f(y)|=|\phi_1(f)+\phi_2(f)|/2\le (|\phi_1(f)|+|\phi_2(f)|)/2\le 1
\]
we have $\phi_1(f)=\phi_2(f)$. By $\phi_0=(\phi_1+\phi_2)/2$ we infer that $f(y)=\phi_0(f)=\phi_1(f)=\phi_2(f)$. Thus $\phi_1,\phi_2\in L$. As $\phi_0\in \ext L$ we have $\phi_0=\phi_1=\phi_2$. Thus $\phi_0\in \ext\ball(E^*)$. 
By Corollary \ref{akc} there exists $x\in Y$ and a unimodular complex number $\lambda$ such that $\phi_0=\lambda \tau_x$ on $E$. Note that $x\in \ch(E)$ since $\bar\lambda \phi_0=\tau_x$ is an extreme point of $\ball(E^*)$ for $\lambda$ is a unimodular complex number. We have 
\[
|f(x)|=|\lambda\tau_x(f)|=|\phi_0(f)|=|f(y)|=\|f\|.
\]
\end{proof}
\begin{prp}\label{5in31jan}
Let $E$ be a $\K$-linear subspace of $C_0(Y,\K)$. 
Suppose that $x\in Y$ is a strong boundary point.
Then the representing measure for $\tau_x$ is only $D_x$
\end{prp}
\begin{proof}
Suppose that $x\in Y$ is a strong boundary point for $E$ and $m$ is its representing measure.
Let $U$ be an open neighborhood of $x$. Since $x$ is a strong boundary point, there is a function $f\in E$ with $f(x)=1=\|f\|_\infty$ and $|f|<1$ on $Y\setminus U$.  Thus $1=|\tau_x(f)|\le \|\tau_x\|\le 1$, so we have $1=\|\tau_x\|=\|m\|$. 

We prove $\suppm{m}=\{x\}$. Suppose contrary that there exists $y\in \suppm(m) \setminus\{x\}$. 
Then there exists an open neighborhood $V$ of $y$ and an open neighborhoof $W$ of $x$ such that  $W\cap V=\emptyset$. There exists $g\in E$ such that $g(x)=1=\|g\|_\infty$, $|g|<1$ on $Y\setminus W$.  
Since $Y\setminus W$ is closed set, there exists $\delta>0$ such that $|g|\le 1-\delta$ on $Y\setminus W$. As $W\cap V=\emptyset$ 
 we have $|g|\le 1-\delta$ on $V$. Since $y\in \suppm(m)$ we have $|m|(V)>0$. Since $\tau_x(g)=\int gdm$, we have 
\[
1=|\tau_x(g)|\le \int_V|g|d|m|+\int_{Y\setminus V}|g|d|m|
\le (1-\delta)|m|(V)+|m|(Y\setminus V)<1,
\]
which is a contradiction proving that $\suppm(m)\setminus \{x\}=\emptyset$. As $m$ is a regular measure, $\suppm(m)$ is not empty, so $\suppm(m)=\{x\}$. Thus $m=\lambda D_x$ for some unimodular complex number $\lambda$. As
\[
1=g(x)=\tau_x(g)=\int gdm=\int gd(\lambda D_x)=\lambda g(x),
\]
we have that $\lambda =1$ and $m=D_x$
\end{proof}
\begin{prp}\label{6in31jan}
Let  $E$ be a $\K$-linear subspace of $C_0(Y,\K)$.
Let $x\in Y$.
Suppose that the representing measure for $\tau_x$ is only $D_x$.
Then $x\in \ch(E)$.
\end{prp}
\begin{proof}
Suppose that $x\in Y$ and the representing measure for $\tau_x$ is only $D_x$.
Let $\tau_x=(\phi_1+\phi_2)/2$ for $\phi_1,\phi_2\in \ball(E^*)$. Suppose that $m_j$ is a representing measurer for $\phi_j$ for $j=1,2$.
 Then $(m_1+m_2)/2$ is a representing measure for $\tau_x$. Thus
\[
|\tau_x(f)|=\left|\int fd(m_1+m_2)/2\right|\le \|f\|_\infty \|(m_1+m_2)/2\|
\]
for every $f\in E$. Thus
\[
1=\|\tau_x\|\le \|(m_1+m_2)/2\|\le (\|m_1\|+ \{m_2\|)/2=1.
\]
Hence $\|(m_1+m_2)/2\|=\|\tau_x\|$, so $(m_1+m_2)/2$ is the representing measure for $\tau_x$.
Thus $D_x=(m_1+m_2)/2$ and 
\[
1=D_x(\{x\})=(m_1(\{x\})+m_2(\{x\}))/2.
\]
As $|m_j(\{x\})|\le 1$ for $j=1,2$, we have that $m_j(\{x\})=1$ for $j=1,2$. As $\|m_j\|=1$ for $j=1,2$ we infer that $m_1=m_2=D_x$ and $\tau_x=\phi_1=\phi_2$. We conclude that $\tau_x\in \ext(\ball (E^*))$, so $x\in \ch(E)$.
\end{proof}

The following corollary is  straightforward from Propositions \ref{5in31jan}, \ref{6in31jan}.
\begin{crl}\label{418}
Suppose that $E$ is a $\K$-linear subspace of $C_0(Y,\K)$. Suppose that $x\in Y$ is a strong boundary point. Then $x\in \ch(E)$.
\end{crl}
The converse of Lemma \ref{6in31jan} does not hold in general.
Let $P(\bar{D})$ be the disk algebra on the closed unit disk $\bar{D}$. Let 
\[
E=\{f\in P(\bar{D}):f(i)=if(1)\}.
\] 
Then $E$ is a uniformly closed $\C$-linear subspace of $C(\bar{D},\C)$ which  separates the points of $\bar{D}$ and has no common zeros. In fact $f(z)=z \in E$ separates the points in $\bar{D}$ and $f(y)=y\ne 0$ for $y\in \bar{D}$ with $y\ne 0$. Put $g(z)=(z-1)(z-i)$. Then $g\in E$ and $g(0)=i\ne 0$. Thus $E$ has no common zeros on $\bar{D}$. Furthermore we have the following.
\begin{prp}
Let $E=\{f\in P(\bar{D}):f(i)=if(1)\}$.
Then  $1\in \ch(E)$, and $D_1$ and $-iD_i$ are representing measures for $\tau_1$. 
\end{prp}
\begin{proof}
Suppose that $\tau_1=(p+q)/2$ for some $p,q\in \ball(E^*)$. As $\tau_1(z)=1$, we infer that $\|\tau_1\|=1$. Since $1=\tau_1(z)=(p(z)+q(z))/2$ and $|p(z)|\le 1$, $|q(z)|\le 1$ we infer that $p(z)=q(z)=1$. Hence $\|p\|=\|q\|=1$. Let $m_p$ be a representing measure for $p$ and $m_q$ a representing measure for $q$. We show that $\suppm(m_p)\subset \{1,i\}$. Suppose not; suppose that there exists $y\in \suppm(m_p)\setminus\{1,i\}$. 

Suppose that $|y|<1$. As $\{1,i\}$ is a peak interpolation set for $P(\bar{D})$ (cf. \cite[p.111]{br} or \cite[Lemma 4.1]{HOST}), there is a function $f\in P(\bar{D})$ such that $f(1)=1,=\|f\|$ and  $f(i)=i$. 
By the maximum absolute value principle for analytic functions, we infer that $|f(y)|<1$. Hence there is a positive integer $n$ such that $|f^{4n+1}(y)|<1/2$. As $if^{4n+1}(1)=f^{4n+1}(i)$ we have $f^{4n+1}\in E$. Then put $h=f^{4n+1}$. 

Suppose that $|y|=1$. As $\{1,i, y\}$ is a peak interpolation set \cite[p.111]{br}, there exists $h\in P(\bar{D})$ such that $f(1)=1=\|f\|$, $h(i)=i$, and $|h(y)|<1/2$. 

Let $U_y$ be an open neighborhood of $y$ such that $|h|<1/2$ on $U_y$. Since $y\in \suppm(m_p)$ and the measure $m_p$ is regular, we have $0<|m_p|(U_y)$.  
Then we get
\[
|p(h)|\le|\int_{U_y}|h|d|m_p|+\int_{\bar{D}\setminus U_y}|h|d|m_p|\le
\frac12 |m_p|(U_y)+|m_p|(\bar{D}\setminus U_y)<1.
\]
Since $|q(h)\le 1$, we get
\[
1=|\tau_1(h)|\le (|p(h)+q(h)|)/2<1,
\]
which is a contradiction proving that $\suppm(m_p)\subset \{1,i\}$. 
Then we infer that there exists two complex numbers $\lambda_p$ and $\mu_p$ with $|\lambda_p|+|\mu_p|=1$ such that $m_p=\lambda_pD_1+\mu_pD_i$ In the same way there exists two complex number $\lambda_q$ and $\mu_q$ with $|\lambda_q|+|\mu_q|=1$ such that $m_q=\lambda_q D_1+\mu_q D_i$. Hence      we obtain that $\frac{\lambda_p+\lambda_q}{2}D_1+\frac{\mu_p+\mu_q}{2}D_i$ is a representing measure for $\frac{p+q}{2}=\tau_1$. Thus
\[
f(1)=\frac{\lambda_p+\lambda_q}{2}f(1)+\frac{\mu_p+\mu_q}{2}f(i)
=
\left(\frac{\lambda_p+\lambda_q}{2}+i\frac{\mu_p+\mu_q}{2}\right)f(1)
\]
for every $f\in E$. Hence
\begin{equation}\label{kani}
\frac{\lambda_p+\lambda_q}{2}+i\frac{\mu_p+\mu_q}{2}=\frac{\lambda_p+i\mu_p}{2}+\frac{\lambda_q+i\mu_q}{2}=1.
\end{equation}
Since 
\[
\left|\frac{\lambda_p+i\mu_p}{2}\right|\le \frac{|\lambda_p|+|\mu_p|}{2}=\frac12
\]
and
\[
\left|\frac{\lambda_q+i\mu_q}{2}\right|\le \frac{|\lambda_q|+|\mu_q|}{2}=\frac12,
\]
we infer from \eqref{kani} that $\lambda_p+i\mu_p=1$. Thus
\[
p(f)=\lambda_pf(1)+\mu_pf(i)=(\lambda_p+i\mu_p)f(1)=f(1)=\tau_1(f)
\]
for every $f\in E$. We have that $\tau_1=p$, so $\tau_1\in \ext(\ball(E^*))$. We conclude that $1\in \ch(E)$. It is evident that $D_1$ and $-iD_i$ are different representing measures for $\tau_1$.
\end{proof}
The strong separation condtion ensures  uniquness of the representing meausre for the point evaluation at a Choquet boundary point. In fact, we have
\begin{prp}\label{7in31jan}
Let $E$ be a $\K$-linear subspace of $C_0(Y,\K)$. Then the following is equivalent.
\begin{itemize}
\item[{\rm(i)}]
$E$ strongly separates the points of $\ch(E)$,
\item[{\rm(ii)}] 
for every $x\in \ch(E)$,
the representing measure for $\tau_x$ is only $D_x$ .
\end{itemize}
\end{prp}
\begin{proof}
We prove (i) implies (ii). 
Suppose that $x\in \ch(E)$ and $m$ a representing measure for $\tau_x$. Note that $\|\tau_x\|=1$ since $\tau_x\in \ext\ball(E^*)$. We have $\|m\|=1$.
Let $y\in \suppm(m)$. We prove that for every open neighborhood $U$ of $y$ we have $|m|(U)=1$. 
(If it were proved, then $m=\lambda D_y$ for some complex number $\lambda$ of modulus $1$ since $m$ is a regular measure. Then 
\begin{equation}\label{sono}
f(x)=\tau_x(f)=\int fdm=\int fd\lambda D_y=\lambda f(y)=\lambda\tau_y(f)
\end{equation}
for every $f\in E$. Then we have $\tau_y=\bar{\lambda}\tau_x$. Since $\tau_x$ is an extreme point of $\ball(E^*)$, $\tau_y$ is also an extreme point of $\ball(E^*)$. Thus $y\in\ch(E)$. By the condition (i)  that $E$ strongly separates the points of $\ch(E)$, we have from \eqref{sono} that $\lambda=1$ and $x=y$. Thus $m=D_x$ follows.)
Suppose not: suppose that there is an open neighborhood $U_0$ of $y$ such that $|m|(U_0)\ne 1$. Then $|m|(U_0)<1$ as $\|m\|=1$. As $m$ is regular and $y\in\suppm(m)$, we have $0< |m|(U_0)$. 
Since $|m|(Y)=\|m\|=1$, $|m|(Y\setminus U_0)|>0$. Put
\[
\varphi_1(f)=\frac{1}{|m|(U_0)}\int_{U_0}f dm,\quad f\in E,
\]
\[
\varphi_2(f)=\frac{1}{|m|(Y\setminus U_0)}\int_{Y\setminus U_0}f dm,\quad f\in E.
\]
It follows that $\varphi_1,\varphi_2\in \ball(E^*)$ and $\tau_x=|m|(U_0)\varphi_1+|m|(Y\setminus U_0)\varphi_2$, where $m|(U_0)>0$, $|m|(Y\setminus U_0)>0$ and  $|m|(U_0)+|m|(Y\setminus U_0)=1$. As $\tau_x\in \ext\ball(E^*)$, we have that $\tau_x=\varphi_1$. In the same way we have 
\[
\tau_x(f)=\frac{1}{|m|(V)}\int_Vfdm, \quad f\in E
\]
for any open neighborhood $V$ of $y$ with $V\subset U_0$. Since $f$ is continuous, for every $\varepsilon>0$, there exists an open neighborhood $V_\varepsilon$ of $y$ such that $|f-f(y)|<\varepsilon$ on $V_\varepsilon$. Hence
\[
|f(y)-f(x)|\le \frac{1}{|m|(V_\varepsilon)}\int_{V_\varepsilon}|f(y)-f|d|m|\le \varepsilon.
\]
Hence $f(y)=f(x)$ for every $f\in E$, so $\tau_y=\tau_x$. As $\tau_x$ is an extreme point of $\ball (E^*)$, so is $\tau_y$. Hence $y\in \ch(E)$. Since $E$ strongly separates the points of $\ch(E)$, we have that $x=y$. It follows that for any $y\in \suppm(m)$, $y$ coincides with $x$, that is, $\suppm(m)=\{x\}$, which is a contradiciton since we assume that $|m|(U_0)<1$. We conclude that $|m|(U)=1$ for any open neighborhood $U$ of $y$. 

We prove (ii) implies (i) by reductio ad absurdum.
Suppose  that $E$ does not strongly separate the points of $\ch(E)$: there exists a pair $x$ and $y$ of different points in $\ch(E)$ such that the equation $|f(x)|=|f(y)|$ holds for every $f\in E$. 
As $\tau_x\in \ext\ball(E^*)$, $\|\tau_x\|=1$ holds, so there exists $f_0\in E$ such that $f_0(x)=1$. Then there exists a complex number $\lambda_0$ with $|\lambda_0|=1$ such that $f_0(y)=\lambda_0 f_0(x)=\lambda_0$. 
For any $f\in E$ with $f(x)\ne 0$ there exists a complex number $\lambda_f$ of unit modulus such that 
$f(y)=\lambda_f f(x)$. As $(f/f(x)+f_0)(x)=2$, we have 
\[
\lambda_f+\lambda_0=(f/f(x)+f_0)(y)=\lambda_{f/f(x)+f_0}(f/f(x)+f_0)(x)=2\lambda_{f/f(x)+f_0}
\]
for every $f\in E$ with $f(x)\ne 0$. As $|\lambda_f|=|\lambda_0|=|\lambda_{f/f(x)+f_0}|=1$ we infer that 
$\lambda_f=\lambda_0$. Thus $f(y)=\lambda_0 f(x)$ for every $f\in E$ with $f(x)\ne 0$. This equation also holds for $f\in E$ with $f(x)=0$. We conclude that $f(y)=\lambda_0 f(x)$ for every $f\in E$. Thus $\bar{\lambda_0}D_y$ is a representing measure for $\tau_x$. As $x\ne y$ we have at least two representing measures $\D_x$ and $\bar{\lambda_0}D_y$ for $\tau_x$. 
\end{proof}
\begin{crl}\label{420}
Let $E$ be a $\K$-linear subspace of $C_0(Y,\K)$. Suppose that $E$ strongly separates the points of $\ch(E)$. Then $x\in Y$ is in the Choquet boundary if and only if the representing measure for $\tau_x$ is only $D_x$.
\end{crl}
\begin{proof}
It is straightforward from Propositions \ref{6in31jan},  \ref{7in31jan}.
\end{proof}
\begin{crl}
Suppose that $E$ is an extremely $C$-regular subspace of $C_0(Y,\K)$. Then $x\in Y$ is in the Choquet boundary if and only if the representing measure for $\tau_x$ is only $D_x$.
\end{crl}
\begin{proof}
By the definition of an extremely $C$-regular subspace, $E$ strongly separates the points of $\ch(E)$. 
Then by Corollary \ref{420} the conclusion holds.
\end{proof}
If $E$ is a subspace of $C(X,\K)$ which contains constants or $E$ is a subalgebra of $C_0(Y,\K)$ which separates the points of $Y$, then the Choquet boundary points are characterized by the uniquness of the  representing measures for the corresponding point evaluations.
\begin{crl}
Supposes that $X$ is a compact Hausdorff space and $E$ is a $\K$-linear subspace of $C(X,\K)$ which separates the points of $X$ and contains constants. 
Then  $x\in Y$ is in the Choquet boundary if and only if the representing measure for $\tau_x$ is only $D_x$. 
\end{crl}
\begin{proof}
  Proopsition \ref{ssp} asserts that $E$ strongly separates the  points of $X$. 
  Then by Corollary \ref{420} we have the conclusion.
\end{proof}
\begin{crl}\label{421}
Suppose that $E$ is a subalgebra of $C_0(Y,\K)$ which separates the points of $Y$. 
Then $x\in Y$ is in the Choquet boundary if and only if the representing measure for $\tau_x$ is only $D_x$. 
\end{crl}
\begin{proof}
 Proopsition \ref{ssp} asserts that $E$ strongly separates the  points of $X$. 
  Then by Corollary \ref{420} we have the conclusion.
\end{proof}
We summarize the results to generalize a theorem of Bishop and de Leeuw  on a characterization of the Choquet boundary for uniform algebras  \cite[p. 39]{phelps}, \cite[Theorem 2.2.6]{br} (cf. \cite[Theorem 2.1]{rrII}, \cite[Theorem 9]{rr}).
\begin{thm}\label{chalg}
Suppose that $A$ is a closed $\K$-subalgebra of $C_0(Y, \K)$ which separates the points of $Y$ and has no common zeros. Let $x\in Y$.The following are equivalent. 
\begin{itemize}
\item[{\rm (i)}]
$x\in \ch(A)$,
\item[{\rm (i$'$)}]
$x\in \ch(\dot{A}+\K)$,
\item[{\rm (ii)}]
$x$ is a strong boundary point for $A$,
\item[{\rm (ii$'$)}]
$x$ is a strong boundary point for $\dot{A}+\K$,
\item[{\rm (iii)}]
the representing measure for the point evaluation $\tau_x$ on $A$ is only $D_x$,
\item[{\rm (iv)}]
there exists a pair of $0<\alpha<\beta\le1$ such that for every open neighborhood $U$ of $x$ there exists a function  $f\in A$ such that $\|f\|_\infty \le 1$, $|f(x)|\ge\beta$, and  $|f|<\alpha$ on $Y\setminus U$,
\item[{\rm (v)}]
for every pair of $0<\alpha<\beta\le1$ and for every open neighborhood $U$ of $x$ there exists a function  $f\in A$ such that $\|f\|_\infty \le 1$, $|f(x)|\ge\beta$, and  $|f|<\alpha$ on $Y\setminus U$.
\end{itemize}
\end{thm}
\begin{proof}
If $\K=\R$, then by the Stone-Weierstrass theorem, $A=C_0(Y,\R)$. It follows that the conditions in  (i) through  (v) hold for every $x\in Y$.

Suppose that $K=\C$.  
(i) $\leftrightarrow$ (iii) is just Corollary \ref{421}.
(ii) $\rightarrow$ (i) is Corollary \ref{418}.
(ii) $\leftrightarrow$ (ii$'$) is Proposition \ref{49}.

 Since $\dot{A}+\C$ is a uniform algebra on $Y_{\infty}$, (i$'$) $\rightarrow$ (ii$'$) follows from the Bishop-de Leeuw theorem (cf. \cite[p.39]{phelps}, \cite[Theorem 2.2.6]{br}).
 
We prove (i) $\rightarrow$ (i$'$). Let $x\in \ch(A)$. To distinguish the point evaluations on $A$ and $\dot{A}+\C$, denote the point evaluation on $\dot{A}+\C$ by $\dot{\tau_x}: \dot{A}+\C\to \C$. The point evaluation on $A$ is denoted $\tau_x$ as usual.  Suppose that $\dot{\tau_x}=(p+q)/2$, where $p,q\in \ball((\dot{A}+\C)^*)$. As $1=\dot{\tau_x}(1)=(p(1)+q(1))/2$ and $|p(1)|\le 1$, $|q(1)|\le 1$ we infer that $p(1)=q(1)=1$. Hence for each $\dot{f}+\lambda\in \dot{A}+\C$ we have
\[
\tau_x(f)+\lambda=\dot{\tau_x}(\dot{f}+\lambda)=(p(\dot{f}+\lambda)+q(\dot{f}+\lambda))/2=
(p(\dot{f})+q(\dot{f}))/2+\lambda.
\]
Define $p':A\to \C$ by $p'(g)=p(\dot{g})$ for $g\in A$, we have $p'\in \ball (A^*)$. In the same way $q'$ is defined  and $q'\in \ball (A^*)$. By the above equality we have $\tau_x=(p'+q')/2$ As $\tau_x\in \ext (\ball(A^*))$, $p'=q'=\tau_x$. It follows that $p=q=\dot{\tau_x}$, proving that $\dot{\tau_x}\in \ext((\dot{A}+\C)^*)$, so $x\in \ch(\dot{A}+\C)$.

We prove (i$'$) $\rightarrow$ (i).
Suppose that $x\in \ch(\dot{A}+\C)$, in other wards, $\dot{\tau_x}\in \ext((\dot{A}+\C)^*)$. 
We have already proved that (i$'$) implies (ii$'$) and  (ii$'$) implies (ii). Hence $x$ is a strong boundary point for $A$. Thus we infer that $1=\|\tau_x\|=|\dot{\tau_x}|\dot{A}\|=1$. 
Let $\Delta_x$ denote a Hahn-Banach extension on $\dot{A}+\C$ of $\dot{\tau_x}|\dot{A}$. 
We show that $\Delta_x=\dot{\tau_x}$. 
By the Riesz-Kakutani theorem there exists 
 $m_x$ of a representing measure of $\Delta_x$ on $Y_\infty$. Note that $\|m_x\|=\|\Delta_x\|=\|\dot{\tau_x}|\dot{A}\|=1$ and $\Delta_x(\dot{f}+\lambda)=\int (\dot{f}+\lambda)dm_x$ for $\dot{f}+\lambda\in \dot{A}+\C$.  Note also that $1=\dot{\tau_x}(1)=\int 1dm_x$ ensures that $m_x$ is a probability measure. 
As $x$ is a strong boundary point for $A$,
by Proposition \ref{4.4} there exists a family $\{K_\alpha\}$ of peak sets for $A$ such that $\bigcap_{\alpha}K_\alpha=\{x\}$. Denote $f_\alpha\in A$ the corresponding peaking function for $K_\alpha$. By the bounded convergence theorem for the probability measure $m_x$ we get
 \[
 1=\Delta_x(f_\alpha^n)=\int f_\alpha^n dm_x \to m_x(K_\alpha)
 \]
 as $n\to 0$ since $f_\alpha=1$ on $K_\alpha$ and $|f_\alpha|<1$ on $Y\setminus K_\alpha$. We see that $m_x(K_\alpha)=1$ for every peak set $K_\alpha$ for $A$ which contains $x$. Let $U$ be an open neighborhood of $x$ in $Y$. Although $Y$ needs not be compact, by considering that $U$ and $K_\alpha's$ are subsets of compact space $Y_\infty$, there exists a finite number of $K_{\alpha_1}, \dots, K_{\alpha_n}$ such that $U\supset \bigcap_{j=1}^nK_{\alpha_j}$. Then we have $1=m_x(\bigcap_{j=1}^nK_{\alpha_j})\le m_x(U)\le 1$. As $U$ is arbitrary open neighborhood of $x$, we get $m_x(\{x\})=1$ since $m_x$ is a regular measure. Thus $m_x=D_x$. Hence $\Delta_x(\dot{f}+\lambda)=\int \dot{f}+\lambda dD_x=\dot{\tau_x}(\dot{f}+\lambda)$ for every $\dot{f}+\lambda\in \dot{A}+\C$. We get that $\Delta_x=\dot{\tau_x}$. 
 
 We prove that $\tau_x\in \ext(\ball(A^*))$. 
 Suppose that $\tau_x=(p+q)/2$ for some $p,q\in \ball (A^*)$. 
 Let $\check{p}:\dot{A}\to \C$ be defined as $\check{p}(\dot{f})=p(f)$, $\dot{f}\in \dot{A}$ and $\check{q}:\dot{A}\to \C$ be defined as $\check{q}(\dot{f})=q(f)$, $\dot{f}\in \dot{A}$. Note that it is well defined since $f\mapsto \dot{f}$ is a bijection from $A$ onto $\dot{A}$. 
  Note also that 
 \[
 \left(\frac{\check{p}+\check{q}}{2}\right)(\dot{f})=\left(\frac{p+q}{2}\right)(f)=\tau_x(f)=\dot{\tau_x}(\dot{f}), \quad f\in A, 
 \]
 so $\dot{\tau_x}|\dot{A}=(\check{p}+\check{q})/2$. 
 Let $\dot{p}$ and $\dot{q}$ be Hahn-Banach extensions of $\check{p}$ and $\check{q}$ on $\dot{A}+\C$ respectively. 
 As $(\dot{p}+\dot{q})/2$ is an extension of $(\check{p}+\check{q})/2=\dot{\tau_x}|\dot{A}$, we have
 \[
 1=\|(p+q)/2\|=\|(\check{p}+\check{q})/2\|\le \|(\dot{p}+\dot{q})/2\|\le 
 (\|\dot{p}\|+\|\dot{q}\|)/2=1.
 \]
 Thus $\|(\dot{p}+\dot{q})/2\|=1$, so $(\dot{p}+\dot{q})/2$ is a Hahn-Banach extension of $\dot{\tau_x}|\dot{A}$. By the result in the previous paragraph that a Hahn-Banach extension of 
 $\dot{\tau_x}|\dot{A}$ is always $\dot{\tau_x}$ we have 
 \[
 \dot{\tau_x}=(\dot{p}+\dot{q})/2.
 \]
 As $\dot{\tau_x}\in \ext(\ball(\dot{A}+\C^*))$ we have $\dot{p}=\dot{q}=\dot{\tau_x}$. 
 For every $f\in A$ we have $\dot{\tau_x}(\dot{f})=\tau_x(f)$, $\dot{p}(\dot{f})=\check{p}(\dot{f})=p(f)$, and $\dot{q}(\dot{f})=\check{q}(\dot{f})=q(f)$ for every $f\in A$, we have $\tau_x(f)=p(f)=q(f)$ for every $f\in A$. We conclude that
 $p=q=\tau_x$, $\tau_x\in \ext(\ball(A^*))$. It follows that  $x\in \ch(A)$.
 
 Suppose that iv) holds. We prove (ii$'$). 
 For every open neighborhood (as a subset of $Y_\infty$) $U$ of $x$ there exists a function  $\dot{f}\in \dot{A}\subset \dot{A}+\C$ such that $\|\dot{f}\|_\infty \le 1$, $|f(x)|\ge\beta$, and  $|f|<\alpha$ on $Y_\infty\setminus U$. Then by \cite[Theorem 2.3.4]{br} we have $x$ is a strong boundary point for $\dot{A}+\C$; (ii$'$) holds. 
 
 Suppose that (ii$'$) holds. We prove (iv). Let $U$ be an open neighborhood of $x$. As we have already pointed out that (ii$'$) is equivalent to (ii), there exists $f\in A$ such that $f(x)=1=\|f\|_\infty$ and $|f|<1$ on $Y\setminus U$. As $Y\setminus U$ is closed, there exsits $\delta<1$ such that $|f|<\delta$ on $Y\setminus U$. Put $\alpha=\delta$ and $\beta=(1+\alpha)/2$. Then we have that $\|f\|=1=f(x)>\beta$ and $|f|<\alpha$ on $Y\setminus U$; iv) holds.

(v) $\rightarrow$ (iv) is trivial. 

We prove (ii) $\rightarrow$ (v).  Suppose that $x$ is a strong boundary point for $A$. Let $\alpha, \beta$ be any pair such that $0<\alpha<\beta\le1$.  Let $U$ be an arbitrary open neighborhood of $x$. Then there exists $f\in A$ such that $f(x)=1=\|f\|$ and $|f|<1$ on $Y\setminus U$. As $Y\setminus U$ is closed, there exists $\delta<1$ such that $|f|<\delta$ on $Y\setminus U$. For a sufficently large positive integer $n$ the inequality $\delta^n<\alpha$ holds. Thus $f^n\in A$ satisfies  $\beta\le|f^n(x)|=1=\|f^n\|_\infty$ and $|f^n|<\alpha$ on $Y\setminus U$.
\end{proof} 
\begin{thm}\label{checr}
Suppose that $E$ is a uniformly closed extremely $C$-regular $\K$-linear subspace of $C_0(Y,\K)$. 
The following are equivalent.
\begin{itemize}
\item[{\rm (i)}]
$x\in \ch(E)$,
\item[{\rm (ii)}]
$x$ is a strong boundary point for $E$,
\item[{\rm (iii)}]
the representing measure for the point evaluation $\tau_x$ on $E$ is only $D_x$,
\end{itemize}
\end{thm}
\begin{proof}
Since $E$ is extremely $C$-regular, $E$ strongly separates the points in $\ch(E)$.
By Corollary \ref{420} we have that (i) $\leftrightarrow$ iii).
The implication (i) $\rightarrow$ (ii) also follows from the definition of the extremely $C$-regularity.
Suppose that (ii) holds. 
 Then by Corollary \ref{418} (i) holds.
\end{proof}
Even if $E$ is a strongly separating $\K$-linear subspace of $C_0(Y,\K)$, a point $y\in \ch(\dot{E}+\K)$ needs not be a point in $\ch(E)$. 
\begin{exa}
Let $E$ be the space defined in Example \ref{0ofE}. Then $1\not\in \ch(E)$ while $1\in \ch(\dot{E}+\K)$. The reason is as follows. 
The space $E$ is strongly separating. Hence $t\in (0,1]$ is in $\ch(E)$ (resp. $\ch(\dot{E}+\K)$) if and only if the representing measure for $\tau_t$ is unique. It is trivial that $D_1$ and $-D_{\frac14}$ are both representing measure for $\tau_1$. Thus $1\not\in \ch(E)$. On the other hand $1$ is a strong boundary point for $\dot{E}+\K$. Thus $1\in \ch(\dot{E}+\K)$.
\end{exa}


\subsubsection{The  \v Silov boundary}\label{sec X}

According to \cite[p.80]{af} the existence of \v Silov boundary for a subalgebra of $C(X,\K)$ which separates the points of $X$ and contains constant is given by \v Silov \cite{silov}.
A simple proof of the existence of the \v Silov boundary for a $\K$-linear subspace of $C(X,\K)$ which separates the points of $X$ and contains constants for a compact Hausdorff  space is exhibited by Bear \cite{bear}. 
In \cite[Proposition 6.4]{phelps} Phelps showed that
 the closure of the Choquet boundary is the \v Silov boundary.  A further references on \v Silov boundary see \cite{af, bl, bk, kaje,rr} for example.
\begin{dfn}\label{shilov}
Let $E$ be a $\K$-linear subspace of $C_0(Y,\K)$. 
We say that a closed subset of $K$ in $Y$ is a {\em \v Silov boundary} for $E$ if it is the smallest closed  boundary in the sense that it is a boundary for $E$ and $K\subset L$ for any {\em closed boundary} $L$ for $E$: $L$ is a boundary for $E$  which is a closed subset of $Y$.
\end{dfn}
Araujo and Font \cite[Theorem 1]{af} proved  that the closure of the Choquet boundary for $E$ is the \v Silov boundary for $E$ if $E$ is a $\K$-linear subspace of $C_0(Y,\K)$ which strongly separates the points of $Y$. The following a little bit generalizes it.
\begin{prp}
Let $E$ be a $\K$-linear subspace of $C_0(Y,\K)$. If $E$ strongly separates the points of $\ch(E)$, then the closure $\overline{\ch(E)}$ of the Choquet boundary in $Y$ is the \v Silov boundary for $E$.
\end{prp}
\begin{proof}
Let $K$ be a closed boundary for $E$. Let $x\in \ch(E)$. We prove that  $x\in K$. 
As $K$ is a boundary, the restriction $E|K$ of $E$ on  $K$ is uniformly closed $\K$-subspace of $C_0(K,\K)$. The restriction map $T:E\to E|K$ by $T|f)=f|K$, $f\in E$ is a bijection and an isometry.
We define $\dot{\tau_x}:E|K\to \K$ by $\dot{\tau_x}(F)=\tau_x(T^{-1}(F))$, $F\in E|K$. Then 
\[
T^*\circ \dot{\tau_x}(f)=\dot{\tau_x}(T(f))=\tau_x(T^{-1}(T(f)))=\tau_x(f),\quad f\in E,
\]
so we infer that $T^*\circ \dot{\tau_x}=\tau_x$ on $E$. 

We prove that $\dot{\tau_x}\in \ext((E|K)^*)$. Suppose that $\dot{\tau_x}=(p+q)/2$ for $p,q\in (E|K)^*$. 
Then
\[
\tau_x=T^*\circ \dot{\tau_x}=\frac{T^*\circ p+T^*\circ q}{2}.
\]
As $\tau_x\in \ext E^*$ we have $T^*\circ p=T^*\circ q=\tau_x.$
Hence $p=q=\dot{\tau_x}$ since $T^*$ is a bijection. It follows that $\dot{\tau_x}\in \ext((E|K)^*)$.

By Corollary \ref{akc} there exist $y\in K$ and $\lambda \in \T$ such that $\dot{\tau_x}=\lambda\tau_y|(E|K)$, where $\tau_y|(E|K):E|K\to \K$ by $(\tau_y|(E|K))(F)=F(y)$ for each $F\in E|K$. It is well defined since $y\in K$. We have
\[
(T^*\circ (\tau_y|(E|K)))(f)=(\tau_y|(E|K))(T(f))=(\tau_y|(E|K))(f|K)=f(y)=\tau_y(f),\quad f\in E,
\]
so
\[
\tau_x(f)=(T^*\circ\dot{\tau_x})(f)=(T^*\circ (\lambda\tau_y)|(E|K)))(f)=\lambda\tau_y(f),\quad f\in E.
\]
Thus $\tau_x=\lambda\tau_y$ on $E$. It follows that $D_x$ and $\lambda D_y$ are representing measures for $\tau_x$. By Proposition \ref{7in31jan} we see that a representing measure for $\tau_x$ is only $D_x$ since we assume $E$ strongly separates the points of $\ch(E)$. We conclude that $x=y$ and $\lambda=1$. Thus $x\in K$, $\ch(E)\subset K$, so the closure $\overline{\ch(E)}$ of $\ch(E)$ is a subset of $K$. As $\ch(E)$ is a boundary for $E$ (Proposition \ref{chbou}) so is $\overline{\ch(E)}$. We conclude  that $\overline{\ch(E)}$ is the \v Silov boundary.
\end{proof}
It is not always the case that the \v Silov boundary exists.
A simple example is as follows.
\begin{exa}
Let $E=\{f\in C(\T,\C): f(\lambda )=\lambda f(1), \,\,\lambda\in \T\}$. Then $E$ is $\C$-linear subspace of $C(\T,\C)$ which separates the points in $\T$. It is easy to see that $\ch(E)=\T$. On the other hand $\{\lambda\}$ is a closed boundary for $E$, for each $\lambda\in \T$. Thus there is no smallest closed boundary for $E$.
\end{exa}
Even if the \v Silov boundary exists, it needs not coincide with the closure of the Choquet boundary.
\begin{exa}
Let $X=[0,1]\cup\{2\}$. Let $E=\{f\in C(X,\K):f(2)=-f(1)\}$. It is evident that $[0,1]$ is the \v Silov boundary and $\ch(E)=X$. Note that $\ch(E)$ separates, but does not strongly separate,  $1$ and $2$. 
Note also that $-D_1$ and $D_2$ are representing measures for $\tau_2$
\end{exa}
\begin{crl}\label{431}
Let  $E$ be a $\K$-linear subspace of $C_0(Y,\K)$.  If $E$ strongly separates the points of $Y$, 
then $\overline{\ch(E)}$ is the \v Silov boundary for $E$. In particular, if $Y$ is compact and $E$  contains constants, or $E$ is a subalgebra, then $\overline{\ch(E)}$ is the \v Silov boundary.
\end{crl}
The case of $1\in E$ is described in \cite[Proposition 6.4]{phelps}.
If $Y$ is compact and $E$ is a strongly separating spaces, then \cite[Proposition 6]{rr} described the above corollary.

We prove Proposition \ref{alcreg}.

\begin{newProof}
Let $x\in \ch(A)$. 
Suppose that $U$ is an open neighborhood of $x$. 
Letting $\beta=1$ and $\varepsilon = \alpha$ for (v) of Theorem \ref{chalg}  we assert that there exists $g\in A$ such that $|g(x)|=1=\|g\|_\infty$ and $|g|<\epsilon$ on $Y\setminus U$. Then $f=\overline{g(x)}g$ is the required function which proves that $A$ is extremely $C$-regular. If $\ch(A)$ is closed, then by the definition of extreme regularity, we have that $A|\ch(A)$ is an extremely regular subspace.

Let $S$ be a \v Silov boundary. We prove that $A|S$ is somewhat regular. 
Let $V$ be a non-empty subset of $S$ and $1>\varepsilon>0$ arbitrary. By Corollary \ref{431} there exists $x_0\in \ch(A)\cap V$. Then there exists $f\in A|S$ such that $f(x_0)=1=\|f_0\|_\infty$ and $|f|\le \varepsilon$. Thus $A|S$ is somewhat regular. \hfill \qed
\end{newProof}
The following is a well known example that shows  the Choquet boundary needs not to be closed even  if the \v Silov boundary exists.
\begin{exa}
Let $A=\{f\in P(\bar{D}):f(0)=f(1)\}|\T$, where $p(\bar{D})$ is the disk algebra on the closed unit disk $\bar{D}$. Then $A$ is a uniform algebra on the unit circle $\T$. Then $\ch(A)=\T\setminus \{1\}$ since every point $\lambda\in \T\setminus \{1\}$ is a peak point with the peaking function $(z+\lambda)/2$ while $1$ is not a peak point by the maximum absolute value principle for analytic functions. Note that the \v Silov boundary is $\T$.
\end{exa}

\section{Sets of representatives}\label{sec4}

\subsection{
Is the homogeneous extension linear?}\label{homoex}
Let $B$ be a real or complex Banach space. 
We denote $S(B)$ the unit sphere $\{a\in B:\|a\|=1\}$ of $B$ and $\ball(B)$ the unit ball $\{a\in B:\|a\|\le 1\}$. 
Any singleton $\{a\}$ of $a\in S(B)$ is convex subset of $S(B)$. Applying Zorn's lemma, there exist a maximal convex subset of $S(B)$ which contains $\{a\}$. Hence $S(B)$ is a union of all maximal convex subsets of $S(B)$. 
We denote 
the set of all maximal convex subsets of $S(B)$ by $\F_B$.
Suppose that the map $T:S(B_1)\to S(B_2)$  is  a surjective isometry with respect to the metric induced by the norm, where
$B_1$ and $B_2$ are both real Banach spaces or both complex Banach spaces.  
The homogeneous extension $\widetilde{T}:B_1\to B_2$ of $T$ is defined as
\begin{equation*}
\widetilde{T}(a)=
\begin{cases}
\|a\|T\left(\frac{a}{\|a\|}\right), & 0\ne a\in B_1 \\
0,& a=0.
\end{cases}
\end{equation*}
By the definition $\widetilde{T}$ is a bijection which satisfies  $\|\widetilde{T}(a)\|=\|a\|$ for every $a\in B_1$. The Tingley's problem asks if $\widetilde{T}$ is real-linear or not.

In \cite{hatori} we introduced {\it a set of representatives} which plays a role in study on complex Mazur-Ulam property. For the convenience of the readers we recall it here.
Suppose 
that $F\in \F_B$ for a real or complex Banach space $B$. It is well known that there exists an extreme point $p$ in the closed unit ball $\ball(B^*)$ of the dual space $B^*$ of $B$ such that $F=p^{-1}(1)\cap S(B)$ (cf. \cite[Lemma 3.3]{tanaka}, \cite[Lemma 3.1]{HOST}). Recall that  $\ex (\ball(B^*))$ denotes the set of all extreme points of $\ball(B^*)$.
Put
\[
Q=\{q\in \ex(\ball(B^*)) : q^{-1}(1)\cap S(B)\in \F_B\}.
\]
We define an equivalence relation $\sim$ in $Q$.
We write $\T=\{z\in {\mathbb{C}}:|z|=1\}$ if $B$ is a complex Banach space, where $\C$ denotes the space of all complex numbers, and $\T=\{\pm 1\}$ if $B$ is a real Banach space.

In Definition \ref{def1} through Definition \ref{def5} $B$ is a real or complex Banach space.
\begin{dfn}[Definition 2.1 in \cite{hatori}]\label{def1}
Let  $p_1,p_2\in Q$. We denote $p_1\sim p_2$ if there exits $\gamma \in \T$ such that $p_1^{-1}(1)\cap S(B)=(\gamma p_2)^{-1}(1)\cap S(B)$.
\end{dfn}
Note that $\gamma p\in Q$ if $\gamma\in \T$ and $p\in Q$.
It is a routine argument to show that the binary relation $\sim$ is an equivalence relation in $Q$.
\begin{dfn}[Definition 2.3 in \cite{hatori}]\label{def3}
A set of all representatives with respect to the equivalence relation $\sim$ is 
 simply called a set of representatives for $\F_B$.
\end{dfn}
Note that a set of representatives exists due to the choice axiom. Note also that a set of representatives $P$ for $\F_B$ is a norming family for $B$ in the sense that $\|a\|=\sup_{p\in P}|p(a)|$. Hence it is a uniqueness set for $B$.

\begin{lmm}[Lemma 2.5 in \cite{hatori}]\label{lemma4}
Let $P$ be a set of representatives. For $F\in \F_B$  there exists a unique $(p,\lambda)\in P\times \T$ such that $F=\{a\in S(B):p(a)=\lambda\}$. Conversely, for $(p,\lambda)\in P\times \T$ we have $\{a\in S(B): p(a)=\lambda\}$ is in $\F_B$.
\end{lmm}
For each set of representatives $P$, Lemma \ref{lemma4} gives a bijective  correspondance between $\F_B$ and $P\times \T$.
\begin{dfn}[Definition 2.6 in \cite{hatori}]\label{def5}
For $(q,\lambda)\in  Q\times \T$, we denote $F_{q,\lambda}=\{a\in S(B):q(a)=\lambda\}$. A map 
\[
I_B:\F_B\to P\times \T
\]
is defined by $I_B(F)=(p,\lambda)$ for $F=F_{p,\lambda}\in \F_B$.
\end{dfn}
By Lemma \ref{lemma4} the map $I_B$ is well defined and bijective. An important theorem of Cheng, Dong and Tanaka states that a surjective isometry between the unit spheres of Banach spaces preserves  maximal convex subsets of the unit spheres. This was first exhibited by Cheng and Dong in \cite[Lemma 5.1]{chengdong} and a complete proof was given by Tanaka \cite[Lemma 3.5]{tanaka2014b}.

In the following  $T:S(B_1)\to S(B_2)$ is a surjective isometry between both real Banach spaces or complex Banach spaces $B_1$ and $B_2$. We denote $P_j$ is a set of representatives for $\F_{B_j}$ for $j=1,2$. 
Applying the theorem of Cheng, Dong and Tanaka, a bijection $\TT:\F_{B_1}\to \F_{B_2}$ is well defined.
\begin{dfn}[Definition 2.7 in \cite{hatori}]\label{def7}
The map $\TT:\F_{B_1}\to \F_{B_2}$ is defined by $\TT(F)=T(F)$ for $F\in \F_{B_1}$. 
The map $\TT$ is well defined and bijective. Put 
\[
\Psi=I_{B_2}\circ \TT\circ I_{B_1}^{-1}:P_1\times \T \to P_2\times \T.
\]
Define two maps
\[
\phi:P_1\times \T \to P_2
\]
and
\[
\tau:P_1\times \T \to \T
\]
by 
\begin{equation}\label{psi}
\Psi(p,\lambda)=(\phi(p,\lambda), \tau(p,\lambda)), \,\,(p,\lambda)\in P_1\times \T.
\end{equation}
If $\phi(p,\lambda)=\phi(p,\lambda')$ for every $p\in P_1$ and $\lambda, \lambda'\in \T$ we simply write $\phi(p)$ instead of $\phi(p,\lambda)$.
\end{dfn}
An equivalent form of  \eqref{psi} is as follows.
\begin{equation}\label{fundamental}
T(F_{p,\lambda})=F_{\phi(p,\lambda),\tau(p,\lambda)},\quad (p,\lambda)\in P_1\times \T.
\end{equation}
Note that 
\begin{equation}\label{--}
\phi(p,-\lambda)=\phi(p,\lambda), \,\,
\tau(p,-\lambda)=-\tau(p,\lambda)
\end{equation}
for every $(p,\lambda)\in P_1\times \T$ (cf. \cite{hatori}).
The reason is as follows. First it is well known that $T(-F)=-T(F)$ for every $F\in \F_{B_1}$ (cf. \cite[Proposition 2.3]{mori}). Hence
\begin{multline*}
F_{\phi(p,-\lambda),\tau(p,-\lambda)}=T(F_{p,-\lambda})=T(-F_{p,\lambda})
=-T(F_{p,\lambda}) 
=-F_{\phi(p,\lambda),\tau(p,\lambda)}=F_{\phi(p,\lambda),-\tau(p,\lambda)}
\end{multline*}
for every $p\in P_1$ since $F_{p,-\lambda}=-F_{p,\lambda}$ by the definition of $F_{p,\lambda}$. Since the map $I_{B_2}$ is a bijection we have \eqref{--}.

Rewriting \eqref{fundamental} 
we get  an essential equation in our argument :
\begin{equation}\label{important}
\phi(p,\lambda)(T(a))=\tau(p,\lambda),\quad a\in F_{p,\lambda}.
\end{equation}
Under certain additional assumptions we have
\[
\phi(p,\lambda)=\phi(p,\lambda'), \quad p\in P_1
\]
for every $\lambda$ and $\lambda'$ in $\T$, and 
\[
\tau(p,\lambda)=\tau(p,1)\times
\begin{cases}
\text{$\lambda$, \quad for some $p\in P_1$}
\\
\text{$\overline{\lambda}$, \quad for other $p$'s}
\end{cases}
\]
for $\lambda\in \T$.
We get, under some condition on $B_1$, that
\begin{equation}\label{then}
\phi(p)(T(a))=\tau(p,1)\times
\begin{cases}
\text{$p(a)$, \quad for some $p\in P_1$}
\\
\text{$\overline{p(a)}$, \quad for other $p$'s}
\end{cases}
\end{equation}
for $a\in F_{p,p(a)}$. 
We empharsize that \eqref{then} holds for $a\in S(B_1)$ with $|p(a)|=1$. 
If the equation \eqref{then} holds for any $a\in S(B_1)$, without the restriction that $a\in F_{p,p(a)}$, then applying the definition of $\widetilde{T}$ we get
\begin{equation}\label{thenthen}
\phi(p)(\widetilde{T}(a))=\tau(p,1)\times
\begin{cases}
\text{$p(a)$, \quad for some $p\in P_1$}
\\
\text{$\overline{p(a)}$, \quad for other $p$'s}
\end{cases}
\end{equation}
for {\em every} $a\in B_1$, with which we infer that
\[
\phi(p)(\widetilde{T}(a+rb))=\phi(p)(\widetilde{T}(a))+\phi(p)(r\widetilde{T}(b))
\]
for every pair $a,b\in B_1$ and every real number $r$. By a further consideration,  we will conclude that $\widetilde{T}$ is real-linear. 
It means that we will arrive at the final positive solution for Tingley's problem if \eqref{then} holds for all $a\in S(B_1)$.


\subsection{The Hausdorff distance condition.}
Recall that the Hausdorff distance $d_H(K,L)$ between non-empty closed subsets $K$ and $L$ of a metric space with  metric $d(\cdot, \cdot)$ is defined by
\[
d_H(K,L)=\max\{\sup_{a\in K}d(a,L), \sup_{b\in L}d(b,K)\}.
\]

\begin{dfn}[Definition 3.2 in \cite{hatori}]
Let $B$ be a complex Banach space and $P$ a set of representatives for $\F_B$. 
We say that $B$ satisfies the  Hausdorff distance condition if the equality
\[
d_H(F_{p,\lambda}, F_{p',\lambda'})=2
\]
for every pair $(p,\lambda)$ and $(p',\lambda')$ in $P\times \T$ such that $p\ne p'$.
\end{dfn}
By Lemma 3.1 in \cite{hatori}, $d_H(F_{p,\lambda}, F_{p'\lambda'})=2$ provided that $F_{p,\lambda}\cap F_{p',\lambda'}\ne \emptyset$.
We can formulate the notion of the condition of the Hausdorff distance in terms of $Q$.
\begin{lmm}
A Banach space $B$ satisfies 
the  Hausdorff distance condition if and only if 
$d_H(F_{q,\lambda},F_{q',\lambda'})=2$ for every pair $q$ and $q'$ of $A$ with $q \not\sim q'$.
\end{lmm}
A proof is a routine argument and is omitted.
\begin{lmm}[Lemma 3.4 in \cite{hatori}]\label{lemma9}
Let $B_j$ be a complex Banach space for $j=1,2$ and
$T:S(B_1)\to S(B_2)$ a surjective isometry.
Suppose that $B_1$ satisfies the  Hausdorff distance condition. Let $P_1$ be a set of representatives for $\F_{B_1}$.
Then we have $\phi(p,\lambda)=\phi(p,\lambda')$ for every $p\in P_1$ and $\lambda, \lambda'\in \T$. Put 
\[
P_1^+=\{p\in P_1:\tau(p,i)=i\tau(p,1)\}
\]
and
\[
P_1^-=\{p\in P_1:\tau(p,i)=\bar{i}\tau(p,1)\}.
\]
Then $P_1^+$ and $P_1^-$ are possibly empty disjoint subsets of  
$P_1$ such that $P_1^+\cup P_1^-=P_1$. Furthermore we have
\[
\tau(p,\lambda)=\lambda \tau(p,1),\quad p\in P_1^+,\lambda\in \T
\]
and
\[
\tau(p,\lambda)=\bar{\lambda}\tau(p,1),\quad p\in P_1^-, \lambda\in \T.
\]
\end{lmm}
\begin{proof}
See the proof of \cite[Lemma 3.4]{hatori}
\end{proof}

\subsection{The set $M_{p,\alpha}$ and the Mazur-Ulam property}

 We exhibit the definition of $M_{p,\alpha}$ for a real or complex Banach space. The case of a complex Banach space is in \cite[Definition 4.1]{hatori}. We denote $\bar\D=\{z\in {\mathbb{K}}:|z|\le 1\}$, where $\mathbb{K}=\mathbb{R}$ if the corresponding Banach space is a real one and 
$\mathbb{K}=\mathbb{C}$ if the corresponding Banach space is 
a complex one.
\begin{dfn}\label{def12}
Let $B$ be a real or complex Banach space and $P$ a set of representatives for $\F_B$. 
For $p\in P$ and $\alpha\in \bar\D$ we denote
\[
M_{p,\alpha}=\{a\in S(B): d(a,F_{p,\alpha/|\alpha|})\le 1-|\alpha|, 
d(a, F_{p,-\alpha/|\alpha|})\le 1+|\alpha|\},
\]
where we read $\alpha/|\alpha|=1$ if $\alpha=0$. 
\end{dfn}
\begin{lmm}[cf. Lemma 4.2 in \cite{hatori}]\label{mpalpha}
Suppose that  $B_j$ is a real or complex  Banach space for $j=1,2$, and $T:S(B_1)\to S(B_2)$ is a surjective isometry. 

If $B_j$ is a real Banach space for $j=1,2$, then we have
\[
T(M_{p,\alpha})=\tau(p,1)M_{\phi(p),\alpha}
\]
for every $(p,\alpha)\in P_1\times \T$.

 If $B_j$ is a  complex Banach space $j=1,2$ and $B_1$ satisfies the Hausdorff distance condition, then we have
\[
T(M_{p,\alpha})=
\begin{cases}
\tau(p,1)M_{\phi(p),\alpha},\quad p\in P_1^+
\\
\tau(p,1)M_{\phi(p),\overline{\alpha}},\quad p\in P_1^-
\end{cases}
\]
for every $(p,\alpha)\in P_1\times \T$. Here $P_1^+$ and $P_1^=$ are dfined as in Lemma \ref{lemma9}.
\end{lmm}
\begin{proof}
According to the definition of the map $\Psi$ we have
\[
T(F_{p,\frac{\alpha}{|\alpha|}})=F_{\phi(p,\frac{\alpha}{|\alpha|}), \tau(p,\frac{\alpha}{|\alpha|})}
\]
and
\[
T(F_{p,-\frac{\alpha}{|\alpha|}})=F_{\phi(p,-\frac{\alpha}{|\alpha|}), \tau(p,-\frac{\alpha}{|\alpha|})}
\]

Suppose that  $B_j$ is a real Banach space  Then by the definition $\T=\{\pm1\}$.
By \eqref{--} we have $\phi(p,1)=\phi(p,-1)$ for every $p\in P_1$. 
Hence $\phi(p,\lambda)$ does not depend on the second term for a real Banach space. We also have $\tau(p,-1)=-\tau(p,1)$ for every $p\in P_1$ by \eqref{--}. It follows that
\[
T(F_{p,\frac{\alpha}{|\alpha|}})=F_{\phi(p),\frac{\alpha}{|\alpha|}\tau(p,1)}=\tau(p,1)F_{\phi(p),\frac{\alpha}{|\alpha|}}
\]
and
\[
T(F_{p,-\frac{\alpha}{|\alpha|}})=F_{\phi(p),-\frac{\alpha}{|\alpha|}\tau(p,1)}=\tau(p,1)F_{\phi(p),-\frac{\alpha}{|\alpha|}}
\]
As $T$ is a surjective isometry we have
\begin{multline*}
d(a,F_{p,\frac{\alpha}{|\alpha|}})=
d(T(a), F_{\phi(p,\frac{\alpha}{|\alpha|}), \tau(p,\frac{\alpha}{|\alpha|})})=
d(T(a),F_{\phi(p),\frac{\alpha}{|\alpha|}\tau(p,1)})
\\
=d(T(a),\tau(p,1)F_{\phi(p),\frac{\alpha}{|\alpha|}})=
d(\tau(p,1)T(a), F_{\phi(p),\frac{\alpha}{|\alpha|}})
\end{multline*}
and
\begin{multline*}
d(a,F_{p,-\frac{\alpha}{|\alpha|}})=
d(T(a),F_{\phi(p,-\frac{\alpha}{|\alpha|}), \tau(p,-\frac{\alpha}{|\alpha|})})
=d(T(a),F_{\phi(p),-\frac{\alpha}{|\alpha|}\tau(p,1)})
\\
=d(T(a), \tau(p,1)F_{\phi(p),-\frac{\alpha}{|\alpha|}})=
d(\tau(p,1)T(a), F_{\phi(p),-\frac{\alpha}{|\alpha|}}).
\end{multline*}
As  $T$ is a bijection we conclude that
\[
\tau(p,1)T(M_{p,\alpha})=M_{\phi(p),\alpha}
\]
for every $p\in P_1$ and $\alpha\in \bar\D$, so
\[
T(M_{p,\alpha})=\tau(p,1)M_{\phi(p),\alpha}
\]
for every $p\in P_1$ and $\alpha\in \bar\D$

A proof for the case where  $B_j$ is a complex Banach space is in \cite[Proof of Lemma 4.2]{hatori}.
\end{proof}

\subsubsection{A sufficient condition for the Mazur-Ulam property : the case of a real Banach space.}

\begin{prp}\label{prop2.2}
Let $B$ be a real Banach space and 
$P$ a set of representatives for $\F_{B}$.
Suppose that
\begin{equation}\label{mpar}
M_{p,\alpha}=\{a\in S(B):p(a)=\alpha\}
\end{equation}
for every $p\in P$ and $-1\le \alpha \le 1$. 
Then $B$ has the Mazur-Ulam property.
\end{prp}
\begin{proof}
Let $B_2$ be a real Banach space and $T:S(B_1)\to S(B_2)$ a surjective isometry.
We first prove the following equation \eqref{rthen} for every $p\in P_1$ and $a\in S(B_1)$ 
without assuming that $|p(a)|=1$.
\begin{equation}\label{rthen}
\phi(p)(T(a))=\tau(p,1)p(a)
\end{equation}
for every $p\in P_1$ and $a\in S(B_1)$ with $|p(a)|=1$.
Let $p\in P_1$ and $a\in S(B_1)$. 
Put $\alpha=p(a)$. Then by \eqref{mpar} $a\in M_{p,\alpha}$. 
We have by Lemma \ref{mpalpha} that 
\[
\phi(p)(T(a))=\alpha\tau(p,1)=\tau(p,1)p(a).
\]
It follows that for the homogeneous extension $\widetilde{T}$ of $T$ we have
\[
\phi(p)(\widetilde{T}(c))=\phi(p)\left(\|c\|T\left(\frac{c}{\|c\|}\right)\right)
=\|c\|\tau(p,1)p\left(\frac{c}{\|c\|}\right)=\tau(p,1)p(c)
\]
for every $0\ne c\in B_1$. As the equality $\phi(p)(\widetilde{T}(0))=\tau(p,1)p(0)$ holds, we obtain for $a,b\in B_1$ and a real number $r$ that
\[
\phi(p)(\widetilde{T}(a+rb)=\tau(p,1)p(a+rb)=\tau(p,1)p(a)+r\tau(p,1)p(b)
\]
and
\[
\phi(p)(\widetilde{T}(a)+r\widetilde{T}(b))=\phi(p)(\widetilde{T}(a))+r\phi(p)(\widetilde{T}(b))
=\tau(p,1)p(a)+r\tau(p,1)p(b).
\]
It follows that
\[
\phi(p)(\widetilde{T}(a+rb))=\phi(p)(\widetilde{T}(a)+r\widetilde{T}(b))
\]
for every $p\in P_1$, $a,b\in B_1$, and every real number $r$. As $\phi(P_1)=P_2$ is a norming family we see that $\widetilde{T}$ is real-linear on $B_1$. As the homogeneous extension is a norm-preserving bijection
  as is described in the subsection \ref{homoex} we complete the proof.
\end{proof}


\subsubsection{A sufficient condition for the complex Mazur-Ulam property : the case of a complex Banach space.}

The case of a complex Banach space is exhibited in Proposition 4.4 in \cite{hatori}.
\begin{prp}[Proposition 4.4 in \cite{hatori}]\label{prop15}
Let $B$ be a complex Banach space and $P$ a set of representatives for $\F_{B}$. 
Assume the following two conditions:
\begin{itemize}
\item[\rm{(i)}] 
$B$ satisfies  the Hausdorff distance condition,
\item[\rm{(ii)}]
$M_{p,\alpha}=\{a\in S(B):p(a)=\alpha\}$ for every $p\in P$ and $\alpha\in \mathbb{D}$.
\end{itemize}
Then $B$ has the complex Mazur-Ulam property.
\end{prp}


\section{Banach spaces which satisfy the condition $(*)$}\label{sec5}
\begin{dfn}\label{dfn*}
Let $B$ be a real or complex Banach space. We say that $B$ satisfies the condition $(*)$ whenever
there exists a set of representative $P$ for $\F_B$ with  the condition :
  for every $p\in P$,  $\varepsilon>0$, and a closed subset $F$ of $P$ with respect to the relative topology induced by the weak*-topology on $B^*$ such that $p\not\in F$, there exists $a\in S(B)$ such that $p(a)=1$ and $|q(a)|\le \varepsilon$ for all $q\in F$.  
\end{dfn}
\begin{exa}\label{example*}
Suppose that $E$ is a uniformly closed $C$-regular $\K$-linear subspace of $C_0(Y,\K)$.
Put $P=\{\tau_x:x\in \ch(E)\}$. By Theorem \ref{checr} every point in $\ch(E)$ is a strong boundary point. Hence $F_{p,\lambda}$ for any $(p,\lambda)\in P\times\T$ is a maximal convex set. Since $p\not\sim q$ for $p,q\in P$ with $p\ne q$, we have that $P$ is a set of representatives.  The condition $(*)$ holds with $P$. We proved that a closed subalgebra of $C_0(Y,\C)$ which separates the points of $Y$ and has no common zeros is an extremely $C$ regular. Thus such an algebra satisfies the condition $(*)$.
\end{exa}
Throughout the section we assume that $B$ is a real or complex Banach space and $P$ is a set of representatives for $\F_B$ with which the conditon $(*)$ is satisfied.
\begin{lmm}\label{1}
Let $p_1, p_2\in P$ with $p_1\ne p_2$. Let $\mu_1,\mu_2\in \T$. For every $\varepsilon>0$ and an open neighborhood $U$ of $\{p_1, p_2\}$ with respect to the relative topology on $P$ induced by the weak*-toplolgy on $B^*$, there exists $h\in B$ such that $\|h\|\le 1+\varepsilon$,  $p_j(h)=\mu_j$ for $j=1,2$, and $|q(h)|\le \varepsilon$ for every $q\in P\setminus U$.
\end{lmm}
\begin{proof}
We may assume that $\varepsilon\le 1/3$. 
Let $V_j$ be an open neighborhood of $p_j$ for $j=1,2$ in $P$ such that $V_1\cup V_2\subset U$ and $V_1\cap V_2=\emptyset$. Choose any positive real number $\delta$ with $0<\frac{4\delta}{1-\delta^2}<\varepsilon$. By the condition $(*)$ there exists $f_j\in B$ such that $p_j(f_j)=1=\|f_j\|$ and $|q(f_j)|\le \delta$ for every $q\in P\setminus V_j$ for $j=1,2$. As $V_1\cap V_2=\emptyset$ and $p_2\in V_2$ we have $p_2\in P\setminus V_1$, so $|p_2(f_1)|\le \delta$. We also have that $|p_1(f_2)|\le \delta$. Hence we infer that $0<1-|p_1(f_2)p_2(f_1)|$. Put 
\[
h_1=\frac{f_1-p_2(f_1)f_2}{1-p_1(f_2)p_2(f_1)}.
\]
Then we infer that $p_1(h_1)=1$ and $p_2(h_1)=0$. By a simple calculation we have
\[
\|h_1\|\le \frac{\|f_1\|+|p_2(f_1)|\|f_2\|}{1-|p_1(f_2)||p_2(f_1)|}\le 1/(1-\delta).
\]
For $q\in P\setminus V_1$ we have
\[
|q(h_1)|\le \frac{|q(f_1)|+|p_2(f_1)||q(f_2)|}{1-|p_1(f_2)||p_2(f_1)|}\le 2\delta/(1-\delta^2).
\]
In a similar way, we have $p_2(h_2)=1$ , $p_1(h_2)=0$, $\|h_2\|\le 1/(1-\delta)$ and 
$|q(h_2)|\le 2\delta/(1-\delta^2)$ for every $q\in P\setminus U_2$, where 
\[
h_2=\frac{f_2-p_1(f_2)f_1}{1-p_1(f_2)p_2(f_1)}.
\]

Put $h=\mu_1h_1+\mu_2h_2$. Then  $p_j(h)=\mu_j$ for $j=1,2$. 
We prove that $\|h\|\le 1+\varepsilon$. Let $q\in V_1$. Then $q\in P\setminus V_2$. Hence
\[
|q(h)|\le |q(h_1)|+|q(h_2)|\le \|h_1\|+2\delta/(1-\delta^2)\le1/(1-\delta)+2\delta/(1-\delta^2)< 1+\varepsilon.
\]
Let $q\in V_2$. Then we have $q\in P\setminus V_1$, and $|q(h)|<1+\varepsilon$ follows.
For $q\in P\setminus (V_1\cup V_2)$. We infer that 
\[
|q(h)|\le |q(h_1)|+|q(h_2)|\le 4\delta/(1-\delta^2)<\varepsilon.
\]
In particular, we have $|q(h)|<\varepsilon$ for every $q\in P\setminus U$ since $V_1\cup V_2\subset U$.
Since $P$ is a norming family we infer that $\|h\|\le 1+\varepsilon$.
\end{proof}

\begin{prp}\label{63}
Let $p_1,p_2\in P$ with $p_1\ne p_2$, and $\mu_1,\mu_2\in \T$. Then there exists $f\in S(B)$ such that $p_j(f)=\mu_j$ for $j=1,2$.
\end{prp}
\begin{proof}
The idea of proof comes from the proof of the Bishop's $\frac{1}{4}$ $-$ $\frac34$ criterion (cf. \cite[Theorem .2.3.2]{br}).
We define inductively a sequence $\{U_n\}$ of open  (with respect to the relative topology on $P$ induced by the weak*-topology) neighborhoods of 
 $\{p_1, p_2\}$,  and  a sequence $\{h_n\}$ in $B$ as follows : 
let $\varepsilon$ be as $0<\varepsilon\le 1/3$. Let $U_1$ be any open neighborhood of $\{p_1,p_2\}$. Then by Lemma \ref{1} there exists $h_1\in B$ such that $\|h_1\|\le 1+\varepsilon$, $p_l(h_1)=\mu_l$ for $l=1,2$, and $|q(h_1)|\le \varepsilon$ for every $q\in P\setminus U_1$. Having defined $U_1,\cdots, U_{n-1}$ and $h_1,\cdots,h_{n-1}$, set
\[
U_n=\{\{q\in U_{n-1}:|q(h_j)|<1+2^{-n}\varepsilon, 1\le j\le n-1\}.
\]
By Lemma \ref{1} there exists $h_n\in B$ such that $\|h_n\|\le 1+\varepsilon$, $p_l(h_n)=\mu_l$ for $l=1,2$, and $|q(h_n)|\le \varepsilon$ for every $q\in P\setminus U_n$. Now let 
\[
\hh=\sum_{n=1}^\infty\frac{h_n}{2^n}.
\]
In fact the series converges and $\hh\in B$ since $\|h_n\|\le 1+\varepsilon$ for every $n$. We have that  
\[
p_l(\hh)=p_l\left(\sum_{n=1}^\infty\frac{h_n}{2^{n}}\right)=\sum_{n=1}^\infty\frac{p_l(h_n)}{2^n}=\mu_l
\]
for $l=1,2$. 

To prove $\|\hh\|\le 1$, it is enough to observe $|q(\hh)|\le 1$ for every $q\in P$ since $P$ is a norming family. 
We consider three cases: i)  $q\in P\setminus \bigcup_{n=1}^\infty U_n$; ii) there exists $n$ such that $q\in U_n\setminus U_{n+1}$; iii) $q\in \bigcap_{n=1}^\infty U_n$. It is possible since $U_n\supset U_{n+1}$ for every $n$. 

Suppose that i) occurs.
We have $|q(h_n)|\le \varepsilon$ for every $n$, hence $|q(\hh)|\le \varepsilon\le 1/3$. 

Suppose that ii) occurs. If $q\in U_1\setminus U_2$, we have $q\not\in U_m$ for $m\ge 2$ since $\{U_n\}$ is decreasing. Thus $|q(h_1)|\le \|h_1\|\le 1+\varepsilon$ and $|q(h_m)|\le \varepsilon$ for every $m\ge 2$. Therefore we have that 
\[
|q(\hh)|\le \frac{1+\varepsilon}{2}+\sum_{m=2}^\infty\frac{\varepsilon}{2^m}=1/2+\varepsilon<1
\]
since $\varepsilon\le 1/3$. If $q\in U_n\setminus U_{n+1}$ for some $n\ge 2$, then $|q(h_j)|<1+2^{-n}\varepsilon$ for $1\le j\le n-1$. Since $q\not\in U_{n+1}$ we have that $q\not\in U_k$ for $k\ge n+1$. Hence $|q(h_k)|\le\varepsilon$ for every $k\ge n+1$. Therefore we get
\begin{multline*}
|q(\hh)|\le \sum_{j=1}^{n-1}\frac{|q(h_j)|}{2^j}+\frac{|q(h_n)|}{2^n}+\sum_{k=n+1}^\infty\frac{|q(h_k)|}{2^k} \\
\le (1+2^{-n}\varepsilon)(1-2^{-(n-1)})+(1+\varepsilon)2^{-n}+2^{-n}\varepsilon\le 1
\end{multline*}
since $\varepsilon\le 1/3$. 

Suppose that iii) occurs. We have $|q(h_j)|<1+2^{-n}\varepsilon$ for all $n>j$. Hence $|q(h_j)|\le 1$ for all $j$, so $|q(\hh)|\le 1$. We conclude that $|q(\hh)|\le 1$ for every $q\in P$. 

As $P$ is a norming family we infer that $\|\hh\|\le 1$. Thus $\|\hh\|=1$ since $1=|p_1(\hh)|\le \|\hh\|$.
\end{proof}
The following proposition is a version of an additive Bishop's lemma. The proof of one we proved in \cite[Lemma 5.3]{HOST} requires the existence of the constants in the target algebra. Proposition \ref{64} is a generalization of Lemma 5.3 in \cite{HOST} which provide for every Banach space  with the condition $(*)$.

We read $\frac{\alpha}{|\alpha|}=1$ provided that $\alpha=0$.
\begin{prp}\label{64}
Let $p\in P$ and $f\in \ball(B)$. 
For every $0<r<1$ there exist $H$ and $H'$ in $B$ such that
\begin{equation}\label{+}
  \begin{split} %
    \|H+rf\|&=1, \\
p(H+rf)&=\frac{p(f)}{|p(f)|}, \\
\|H\|&\le 1-r|p(f)|,
  \end{split}   %
\end{equation}
and
\begin{equation}\label{-}
\begin{split}
\|H'+rf\|&=1, \\
p(H'+rf&)=-\frac{p(f)}{|p(f)|}, \\
\|H'\|&\le 1+|p(f)|.
\end{split}
\end{equation}
\end{prp}
\begin{proof}
Put $\alpha=p(f)$. 
The proof of of \eqref{+} is  similar to that of Lemma 5.3 in \cite{HOST}, but a small revision is required since  $B$ is not closed under the multiplication.
The proof of \eqref{-} is different from that of Lemma 5.3 in \cite{HOST}. It requires substantial changes. 

We first prove  \eqref{+}. 
Let $\varepsilon$ be any real number such that $0<\varepsilon<1-r|\alpha|$.
Put
\[
F_0=\{q\in P: |r\alpha-rq(f)|\ge\varepsilon/4\}
\]
and 
\[
F_n=\{q\in P:\varepsilon/2^{n+2}\le|r\alpha-rq(f)|\le\varepsilon/2^{n+1}\}
\]
for a positive integer $n$.  By the condition $(*)$, for every positive integer $n$ there exists $u_n\in B$ such that 
\[
p(u_n)=1=\|u_n\|
\]
and
\[
|q(u_n)|\le \min\left\{\frac{1-r}{1-r|\alpha|}, \frac{1}{2^{n+1}}\right\}
\]
for every $q\in F_0\cup F_n$. Put 
\[
u=\sum_{n=1}^\infty\frac{u_n}{2^n}.
\]
Then $u\in \ball(B)$  since $\|u_n\|=1$ for every $n$. Since $p(u)=\sum_{n=1}^\infty\frac{p(u_n)}{2^n}=1$, we see that $\|u\|=1$. 

Letting $H=\left(\frac{\alpha}{|\alpha|}-r\alpha\right)u$, we prove that $\|H+rf\|\le 1$ in the following three cases: (i) $q\in F_0$; (ii) $q\in F_n$ for some $n\ge1$; (iii) $q\in P\setminus \bigcup_{k=0}^\infty F_k$.

(i) Let $q\in F_0$.  We have 
\[
|q(H+rf)|\le \left|\frac{\alpha}{|\alpha|}-r\alpha\right||q(u)|+r|q(f)|\le (1-r|\alpha|)\frac{1-r}{1-r|\alpha|}+r=1.
\]

(ii) Let $q\in F_n$ for some $n\ge1$.  In this case we  have
\[
|q(u)|\le \sum_{m\ne n}\frac{|q(u_m)|}{2^m}+\frac{|q(u_n)|}{2^n}\le 1-\frac{1}{2^n}+\frac{1}{2^n2^{n+1}}.
\]
As $0<\varepsilon<1-r|\alpha|$ we have 
\begin{multline*}
|q(H+rf)|\le (1-r|\alpha|)|q(u)|+r|\alpha|+ |rq(f)-r\alpha|
\\
\le
(1-r|\alpha|)(1-1/2^n+1/(2^n2^{n+1}))+r|\alpha|+\varepsilon/2^{n+1}
\\
<(1-r|\alpha|)(1-1/2^n+1/(2^n2^{n+1})+1/2^{n+1})+r|\alpha|<1
\end{multline*}

(iii) Let $q\in P\setminus \bigcup_{n=0}^\infty F_n$. In this case we have $q(f)=\alpha$, hence
\[
|q(H+rf)|\le (1-r|\alpha|)|q(u)|+r|\alpha|\le 1-r|\alpha|+r|\alpha|=1.
\]
By (i), (ii) and (iii) we have $\|H+rf\|\le 1$ since $P$ is a norming family. 
As we will see that $p(H+rf)=\alpha/|\alpha|$, it will follow that $\|H+rf\|=1$. 

By simple calculations we have
\[
\|H\|=\left\|\left(\frac{\alpha}{|\alpha|}-r\alpha\right)u\right\|=\left|\fraca - r\alpha\right|\|u\|=1-r|\alpha|
\] 
and 
\[
p(H+rf)=p\left(\left(\fraca-r\alpha\right)u+rf\right)=\left(\fraca-r\alpha\right)p(u)+rp(f)=\fraca.
\]
We have completed the proof of \eqref{+}.

Next we prove \eqref{-}. 
 Suppose that $|\alpha|=1$. Put $H'=-(1+r)f$. Then we have $\|H'\|\le 1+r\le 1+|\alpha|$. We also have that $p(H'+rf)=-P(f)=-\fraca$ and $\|H'+rf\|=\|-f\|\le 1$. Since $1=|p(H'+rf)|\le\|H'+rf\|$, we infer that $\|H'+rf\|=1$. Thus \eqref{-} holds if $|\alpha|=1$.
Suppose that $\alpha=0$. As $\|-f\|\le 1$ and $p(-f)=0=\alpha$, we have by \eqref{+} that there exists $H\in B$ such that $\|H-rf\|=1$, $p(H-rf)=\fraca$ and $\|H\|\le 1-r|\alpha|=1$. Letting $H'=-H$ we have that $\|H'+rf\|=1$, $p(H'+rf)=-1=-\fraca$ and $\|H'\|=\|H\|\le 1=1+|\alpha|$. Thus \eqref{-} holds if $\alpha=0$.

We assume  $0<|\alpha|<1$. To prove \eqref{-} we apply induction. 
Define a sequence $\{a_n\}$ by $a_1=1/3$, $a_{n+1}=(a_n+1)/2$ for every positive integer $n$. 
Put $I_1=(0,a_1]$ and  $I_n=(a_{n-1},a_n]$ for each positive integer  $n\ge 2$. 
For each positive integer $n$, put 
\[
C_n=\{g\in \ball(B): |p(g)|\in I_n\}.
\]
By induction on $n$ we prove that 
for any $g\in C_n$ and $0<s<1$, there exists $H_{g,s}\in B$ such that 
\begin{equation}\label{gs}
\begin{split}
\|H_{g,s}+sg\|&=1, \\
p(H_{g,s}+sg&)=-\frac{p(g)}{|p(g)|}, \\
\|H_{g,s}\|&\le 1+|p(g)|.
\end{split}
\end{equation}
If it will be proved, then \eqref{gs} will hold for every $g\in \ball{B}$ such that $0<|p(g)|<1$ since $\lim_{n\to\infty}a_n=1$. Combining with the results for $\alpha=0$ or $|\alpha|=1$, it will follows that \eqref{-} holds for every $f\in \ball(B)$ and $0<r<1$.

Suppose that $g_1\in C_1$  and $0<s_1<1$. Put $\alpha_1=p(g_1)$. Then $0<|\alpha_1|\le a_1$.
Let $\varepsilon_1$ be as $0<\varepsilon_1<1-\sqrt{s_1}$ 
and 
\[
K_1=\{q\in P:|\alpha-q(g_1)|\ge \varepsilon_1/2\}.
\]
By the condition $(*)$ there exists $v_1\in S(B)$ such that $p(v_1)=1=\|v_1\|$ and $|q(v_1)|\le \varepsilon_1/2$ for every $q\in K_1$. Look at $\|-2\sqrt{s_1}\alpha_1v_1+\sqrt{s_1}g_1\|$. If $q\in K_1$, then $|q(v_1)|\le \varepsilon_1/2$. Therefore
\[
|q(-2\sqrt{s_1}\alpha_1 v_1+\sqrt{s_1}g_1)|\le 2\sqrt{s_1}|\alpha_1||q(v_1)|+\sqrt{s_1}|q(g_1)|\le 2\sqrt{s_1}|\alpha_1|\varepsilon_1/2+\sqrt{s_1}<1
\]
since $|\alpha_1|\le 1/3$. If $q\in P\setminus K_1$, then $|\alpha_1-q(g_1)|<\varepsilon_1/2$. Thus 
\[
|q(-2\sqrt{s}\alpha_1 v_1+\sqrt{s}g_1)|\le\sqrt{s}|\alpha_1-2\alpha_1 q(v_1)|+\sqrt{s}|q(g_1)-\alpha_1|
\le3\sqrt{s}|\alpha_1|+\sqrt{s}\varepsilon_1/2<1
\]
since since $|\alpha_1|\le 1/3$ and $\varepsilon_1<1-\sqrt{s_1}$. We conclude that 
 \[
 \|-2\sqrt{s_1}\alpha_1 v_1+\sqrt{s_1}g_1\|\le 1
 \]
  since $P$ is a norming family. Put $g_1'=-2\sqrt{s_1}\alpha_1 v_1+\sqrt{s_1}g_1$ Then $p(g_1')=-\sqrt{s_1}\alpha_1$. 
Applying \eqref{+} with $\sqrt{s_1}$ and $g_1'$ instead of $r$ and $f$ respectively, we find
 $H_+\in B$ such that 
  \[
  \|H_++\sqrt{s_1}g_1')\|=1, \quad p(H_++\sqrt{s_1}g_1')=\frac{p(g_1')}{|p(g_1')|}=-\frac{\alpha_1}{|\alpha_1|}
  \]
  and
  \[
  \|H_+\|\le 1-\sqrt{s_1}|p(g_1')|=1-s_1|\alpha_1|.
  \]
Letting $H_{g_1,s_1}=H_+-2s_1\alpha_1 v_1$ we infer that 
\[
H_{g_1,s_1}+sg_1=
H_+-2s_1\alpha_1v_1+s_1g_1
=H_++\sqrt{s_1}(-2\sqrt{s_1}\alpha_1 v_1+\sqrt{s_1}g_1)=
H_++\sqrt{s_1}g_1'.
\]
 Hence
\[
\|H_{g_1,s_1}+s_1g_1\|=1,\quad p(H_{g_1,s_1}+s_1g_1)=-\frac{\alpha_1}{|\alpha_1|}.
\]
We also have
\[
\|H_{g_1,s_1}\|\le \|H_+\|+2s_1|\alpha_1|\|v_1\|\le 1-s_1|\alpha_1|+2s_1|\alpha_1|=1+s_1|\alpha_1|\le1+|\alpha_1|.
\]
We have proved that \eqref{gs} holds for $n=1$

Suppose that \eqref{gs} holds for every $1\le n \le m$. 
Let $g_{m+1}\in C_{m+1}$ and $0<s_{m+1}<1$ arbitrary.
Put $\alpha_{m+1}=p(g_{m+1})$. 
Put $0<\varepsilon_{m+1}<1-\sqrt{s_{m+1}}$ and 
\[
K_{m+1}=\{q\in P:|\alpha_{m+1}-q(g_{m+1})|\ge\varepsilon_{m+1}/2\}.
\]
By the condition $(*)$ there exists $v_{m+1}\in B$ such that 
\[
p(v_{m+1})=1=\|v_{m+1}\|
\]
and 
\[
|q(v_{m+1})|\le \frac{\varepsilon_{m+1}}{2(a_{m+1}-a_m)}
\]
for every $q\in K_{m+1}$. Put
\[
f_{m+1}=\sqrt{s_{m+1}}\left(\left(a_m\cdot \frac{\alpha_{m+1}}{|\alpha_{m+1}|}-\alpha_{m+1}\right)v_{m+1}+g_{m+1}\right).
\]
Note that $\left|a_m\cdot \frac{\alpha_{m+1}}{|\alpha_{m+1}|} -\alpha_{m+1}\right|=|\alpha_{m+1}|-a_m$ since $a_m<|\alpha_{m+1}|\le a_{m+1}$. 
Suppose that $q\in K_{m+1}$. Then
\begin{multline*}
|q(f_{m+1})|\le \sqrt{s_{m+1}}(|\alpha_{m+1}|-a_m)|q(v_{m+1})|+\sqrt{s_{m+1}}|q(g_{m+1})|
\\
\le
\frac{\sqrt{s_{m+1}}(|\alpha|-a_m)\varepsilon_{m+1}}{2(a_{m+1}-a_m)}+\sqrt{s_{m+1}}\le
\sqrt{s_{m+1}}(\varepsilon_{m+1}/2+1)<1
\end{multline*}
since $|\alpha_{m+1}|\le a_{m+1}$. 
Suppose that $q\in P\setminus K_1$. Then
\begin{multline*}
|q(f_{m+1})|\le \sqrt{s_{m+1}}\left(\left|a_m\cdot\frac{\alpha_{m+1}}{|\alpha_{m+1}|}-\alpha_{m+1}\right||q(v_{m+1})|+|\alpha_{m+1}|+|q(g_{m+1})-\alpha_{m+1}|\right)
\\
\le
\sqrt{s_{m+1}}((|\alpha_{m+1}|-a_m)+|\alpha_{m+1}|+\varepsilon_{m+1}/2)
=
\sqrt{s_{m+1}}(2|\alpha_{m+1}|-a_m+\varepsilon_{m+1}/2)
\\
\le
\sqrt{s_{m+1}}(2a_{m+1}-a_m+\varepsilon_{m+1}/2)=\sqrt{s_{m+1}}(1+\varepsilon_{m+1}/2)<1
\end{multline*}
Therefore we have that $\|f_{m+1}\|\le 1$ since $P$ is a norming family. 
By a calculation we have $p(f_{m+1})=\sqrt{s_{m+1}}a_m\cdot\frac{\alpha_{m+1}}{|\alpha_{m+1}|}$, hence $|p(f_{m+1})|=\sqrt{s_{m+1}}a_m<a_m$. It means that $|p(f_{m+1})|\in I_k$ for some $1\le k \le m$. By the hypothesis of induction there exists $H_k\in B$ such that 
\[
\|H_k+\sqrt{s_{m+1}}f_{m+1}\|=1, \quad p(H_k+\sqrt{s_{m+1}}f_{m+1})=-\frac{p(f_{m+1})}{|p(f_{m+1})|}=-\frac{\alpha_{m+1}}{|\alpha_{m+1}|}
\]
and
\[
\|H_k\|\le 1+|p(f_{m+1})|=1+\sqrt{s_{m+1}}a_m
\]
Put
\[
H_{g_{m+1},s}=H_k+s_{m+1}\left(a_m\cdot\fraca -\alpha\right)v_{m+1}.
\]
Then $H_{g_{m+1},s_{m+1}}+s_{m+1}g_{m+1}=H_k+\sqrt{s_{m+1}}f_{m+1}$ and 
\[
\|H_{g_{m+1},s_{m+1}}+s_{m+1}g_{m+1}\|=1,\quad p(H_{g_{m+1},s_{m+1}}+s_{m+1}g_{m+1})=-\frac{\alpha_{m+1}}{|\alpha_{m+1}|}.
\]
We also have
\begin{multline*}
\|H_{g_{m+1},s_{m+1}}\|\le \|H_k\|+\left\|s_{m+1}\left(a_m\cdot\frac{\alpha_{m+1}}{|\alpha_{m+1}|} -\alpha_{m+1}\right)v_{m+1}\right\|
\\\le
1+\sqrt{s_{m+1}}a_m+s_{m+1}(|\alpha_{m+1}|-a_m)
<1+\sqrt{s_{m+1}}a_m+\sqrt{s_{m+1}}(|\alpha_{m+1}|-a_m)
\\
\le 1+|\alpha_{m+1}|.
\end{multline*}
We conclude by induction that \eqref{gs} holds if $0<|p(g)|<1$.

\end{proof}


\section{Banach spaces which  satisfy the condition $(*)$ and the Mazur-Ulam property}\label{sec6}
\begin{thm}\label{real*}
Let $B$ be a real Banach space which satisfies the condition $(*)$. 
Then $B$ has the Mazur-Ulam property. 
\end{thm}
\begin{proof}
Let $P$ be a set of representative for $\F_B$ with which  the condition $(*)$ is satisfied. 
We prove that \eqref{mpar} of
Proposition \ref{prop2.2}. 
Let $p\in P$ and $-1\le\alpha\le 1$ arbitrary. 
In the same way as  \cite[Lemma 4.3]{hatori} we infer that
\[
M_{p,\alpha}\subset \{a\in S(B):p(a)=\alpha\}.
\]
We prove the inverse inclusion. 
Suppose that $f\in S(B)$ with $p(f)=\alpha$.
Let $0<r<1$ arbitraty. By Proposition \ref{64} 
there exists $H\in B$ such that $H+rf\in F_{p,\frac{\alpha}{|\alpha|}}$ and $\|H\|\le 1-r|\alpha|$. Thus
\[
\|H+rf -f\|\le 1-r|\alpha|+1-r.
\]
We infer that $d(f, F_{p,\frac{\alpha}{|\alpha|}})\le 1-|\alpha|$.
We also have by Proposition \ref{64} that there exists $H'\in B$ such that $H'+rf\in F_{p,-\frac{\alpha}{|\alpha|}}$ and $\|H'\|\le 1+|\alpha|$. Hence
\[
\|H'+rf-f\|\le 1+|\alpha|+1-r.
\]
We infer that $d(f, F_{p, -\frac{\alpha}{|\alpha|}})\le 1+|\alpha|$. 
Thus $f\in M_{p,\alpha}$. Hence $\{f\in S(B):p(f)=\alpha\}\subset M_{p,\alpha}$. 
It follows from Proposition \ref{prop2.2} that $B$ has the Mazur-Ulam property.
\end{proof}
\begin{crl}\label{ecrmu}
Let $E$ be a uniformly closed extremely $C$-regular $\R$ linear subspace of $C_0(Y,\R)$ for a locally compact Hausdorff space $Y$. Then $E$ has the Mazur-Ulam property.  
In particular, $C_0(Y,\R)$ itself and a uniformly closed extremely regular subspace of $C_0(Y,\R)$ has the Mazur-Ulam property.
\end{crl}
\begin{proof}
We prove that $E$ satisfies the condition $(*)$. 
Put $P=\{\tau_x:x\in \ch(E)\}$. By Theorem \ref{checr} every point in $\ch(E)$ is a strong boundary point. Hence $F_{p,\lambda}$ for any $(p,\lambda)\in P\times\{\pm1\}$ is a maximal convex set. Since $p\not\sim q$ for $p,q\in P$ with $p\ne q$, we have that $P$ is a set of representatives. The condition $(*)$ holds with $P$. Then by Theorem \ref{real*} we see that $E$ has the Mazur-Ulam property.
\end{proof}
\begin{thm}\label{71}
Let $B$ be a complex Banach space which satisfies the condition $(*)$.
Then $B$ has the complex Mazur-Ulam property.
\end{thm}
\begin{proof}
Let $P$ be a set of representative for $\F_B$ with which  the condition $(*)$ is satisfied. 
We prove  (i) and (ii) of Proposition \ref{prop15}. 
Let $F_{p,\lambda},F_{p',\lambda'}\in \F_B$ such that $p\ne p'$. 
By Proposition \ref{63} there exists $f\in S(B)$ such that $p(f)=-\lambda$ and $p'(f)=\lambda'$. Then $f\in F_{p',\lambda'}$. For any $g\in F_{p,\lambda}$ we infer that
\[
2=|p(f)-p(g)|\le \|f-g\|\le  2,
\]
hence $d(f,F_{p,\lambda})=2$. Hence $d_H(F_{p,\lambda}, F_{p',\lambda'})=2$. As a pair $F_{p,\lambda}$ and $F_{p',\lambda'}$ with $p\ne p'$ is arbitrary, we see that $B$ satisfies  the Hausdorff distance condition holds; (i) of Proposition \ref{prop15} holds. 

We prove
 (ii) of Proposition \ref{prop15}.
 Let $p\in P$ and $\alpha\in \bar{D}$ arbitrary. 
Let $0<r<1$ be arbitrary. 
By Proposition \ref{64} 
there exists $H\in B$ with $H+rf\in F_{p,\frac{\alpha}{|\alpha|}}$ and $\|H\|\le 1-r|\alpha|$. There also exists
 $H'\in B$ with  $H'+rf\in F_{p,-\frac{\alpha}{|\alpha|}}$ and $\|H'\|\le 1+|\alpha|$. Hence 
\[
\|H+rf -f\|\le 1-r|\alpha|+1-r,
\]
and 
\[
\|H'+rf-f\|\le 1+|\alpha|+1-r.
\]
As $r$ is arbitrary we infer that $d(f, F_{p,\frac{\alpha}{|\alpha|}})\le 1-|\alpha|$ 
and 
$d(f, F_{p, -\frac{\alpha}{|\alpha|}})\le 1+|\alpha|$.
It follows that  $f\in M_{p,\alpha}$. Thus $\{f\in S(B):p(f)=\alpha\}\subset M_{p,\alpha}$.
The inverse inclusion 
\[
M_{p,\alpha}\subset \{f\in S(B):p(f)=\alpha\}
\]
is 
 by \cite[Lemma 4.3]{hatori}. 
We conclude that (ii) of Proposition \ref{prop15} holds.

It follows from Proposition \ref{prop15} that $B$ has the complex Mazur-Ulam property.
\end{proof}
\begin{crl}\label{ecregmu}
A uniformly closed extremely $C$-regular $\C$-linear subspace of $C_0(Y,\C)$ for a locally compact Hausdorff space has the complex Mazur-Ulam property. 
In particular, $C_0(Y,\C)$ itself and a uniformly closed extremely regular $\C$-linear subspace of $C_0(Y,\C)$ has the complex Mazur-Ulam property.
\end{crl}
\begin{proof}
The proof that $E$ satisfies the condition $(*)$ is essentially the same as that for Corollary \ref{ecrmu} (cf. Example \ref{example*}).  Then by Theorem \ref{71} we see that $E$ has the complex Mazur-Ulam property.
\end{proof}
Hatori \cite[Theorem 4.5]{hatori} proved that a uniform algebra has the complex Mazur-Ulam property. 
In the proof of Theorem 4.5 in \cite{hatori}, it is crucial that a uniform algebra contains the constants. 
Cabezas, Cueto-Avellaneda, Hirota, Miura and Peralta \cite[Corollary 3.2]{cchmp} proved that 
$C_0(Y,\C)$ satisfies the complex Mazur-Ulam property. 
Cueto-Avellaneda, Hirota, Miura and Peralta \cite[Theorem 2.1]{chmp} have proved that each surjective
isometry between the unit spheres of two uniformly closed  algebras on locally compact
Hausdorff spaces which separates the points without common zeros admits an extension to a surjective real linear isometry between these algebras.
The following generalizes a both theorems of Cueto-Avellaneda, Hirota, Miura and Peralta \cite[Theorem 2.1]{chmp} and  Cabezas, Cueto-Avellaneda, Hirota, Miura and Peralta \cite[Corollary 3.2]{cchmp}.
\begin{crl}\label{csacmu}
A (non-zero) closed subalgebra of  $C_0(Y,\C)$ has the complex Mazur-Ulam property.
In particular a uniform algebra has the complex Mazur-Ulam property.
\end{crl}
\begin{proof}
Let $A$ be a non-zero closed subalgebra of $C_0(Y,\C)$. Let 
\[
C=\{y\in Y:\text{$f(y)=0$ for all $f\in A$}\}.
\] 
As we assume $A\ne \{0\}$, there is a non-zero function in $A$, so $C$ is a proper subset of $Y$. For any pair of points  $x$ and $y$ in $Y\setminus C$, we denote  $x\sim y$ if $f(x)=f(y)$ for all $f\in A$. Then $\sim$ is an equivalence relation on $Y\setminus C$. Let $Y_0$ be the quotient space induced by the relation $\sim$. Then $Y_0$ is a  locally compact, possibly compact, Hausdorff  space induced by the quotient topology. We may suppose that $A$ is a closed subalgebra of $C_0(Y_0,\C)$ which separates the points in $Y_0$ and has no common zeros.
Then by Proposition \ref{alcreg} $A$ is a uniformly closed  extremely $C$-regular $\C$-linear subspace of $C_0(Y_0,\C)$. It follows by Corollary \ref{ecregmu}  that $A$ has the complex Mazur-Ulam property.

Suppose that $A$ is a uniform algebra on a compact Hausdorff space. Then $A$ is a closed subalgebra of $C(X,\C)$. Then by the first part we have that $A$ has the complex Mazur-Ulam property.
\end{proof}

\section{Final remarks}\label{sec7}
In section \ref{sec5} we introduced the conidition $(*)$ for Banach spaces and 
 we proved that an extremely  $C$-regular closed subspace of $C_0(Y,\K)$ satisfies the condition $(*)$.  On the other hand we have not enough examples of Banach spaces which satisfy the conditon. 
Cabezas, Cueto-Avellaneda, Hirota, Miura and Peralta \cite{cchmp} preved that every {\em commutative} JB$^*$ triple satisfies the complex Mazur-Ulam property. We do not know if a commutative JB$^*$ triple satisfies the condition $(*)$ or not. If it would satisfy the condition $(*)$, we would get an alternative proof of a theorem of  Cabezas, Cueto-Avellaneda, Hirota, Miura and Peralta \cite[Theorem 3.1]{cchmp}. It is interesting to exhibit examples of Banach spaces which satisfy the condition $(*)$.
We have proved the {\em complex} Mazur-Ulam property especially for a closed subalgebra of $C_0(Y,\K)$ which separates the points of $Y$ and has no common zeros, we expect it also has the Mazur-Ulam property.

\subsubsection*{Acknowledgments.}
I would like to express my thanks to Professor Antonio Peralta for fruitful discussions which started on the end of 2021. It is certain that the discussions and his useful comments are the starting point of the research in this paper. 

Special thanks are due to Professor Shiho Oi for organaizing
the RIMS workshop ``Research on preserver problems on Banach
algebras and related topics'' and for editting RIMS K\^{o}ky\^{u}roku Bessatsu.

This work was supported by JSPS KAKENHI Grant Numbers JP19K03536 and the Research Institute for Mathematical
Sciences, an International Joint Usage/Research Center located in Kyoto
University



\end{document}